\documentclass[a4paper]{article}
\usepackage{kusochek+}
\hypersetup{
pdftitle={Curvature tensor of smoothable Alexandrov spaces},
pdfauthor={Nina Lebedeva and Anton Petrunin}
}

\begin{document}

\title{Curvature tensor of smoothable Alexandrov spaces}
\date{}
\author{Nina Lebedeva and Anton Petrunin} 
\maketitle

\begin{abstract}
We prove weak convergence of curvature tensors of Riemannian manifolds 
for converging noncollapsing sequences with a lower bound on sectional curvature.
\end{abstract}

\newcommand{\Addresses}{{\bigskip\footnotesize

\noindent Nina Lebedeva,
\par\nopagebreak
 \textsc{St. Petersburg State University, 7/9 Universitetskaya nab., St. Petersburg, 199034, Russia}
\par
\nopagebreak
 \textsc{St. Petersburg Department of V. A. Steklov Institute of Mathematics of the Russian Academy of Sciences, 27 Fontanka nab., St. Petersburg, 191023, Russia}
  \par\nopagebreak
  \textit{Email}: \texttt{lebed@pdmi.ras.ru}

\medskip

\noindent   Anton Petrunin, 
\par\nopagebreak
 \textsc{Math. Dept. PSU, University Park, PA 16802, USA.}
  \par\nopagebreak
  \textit{Email}: \texttt{petrunin@math.psu.edu}
  
}}

\section{Introduction}

The weak convergence and measure-valued tensor used in the following theorem are defined in the next section;
a more precise formulation is given in \ref{main}.

\begin{thm}{Main theorem}
Let $M_1,M_2,\dots$ be a sequence of complete $m$-dimensional Riemannian manifolds with sectional curvature bounded below by~$\kappa$.
Assume that the sequence $M_n$ Gromov--Hausdrorff converges to an Alexandrov space $A$ of the same dimension.
Then the curvature tensors of $M_n$ weakly converge to a measured-valued tensor on $A$.
\end{thm}

Note that from the theorem we get that the limit tensor of the sequence depends only on $A$ and does not depend on the choice of the sequence $M_n$.
Indeed, suppose another sequence $M_n'$ satisfies the assumptions of the theorem.
If the limit tensor is different, 
then a contradiction would occur for the alternated sequence $M_1,M_1',M_2,M_2',\dots$ 
In~particular, if the limit space is Riemannian, then the limit curvature tensor is the curvature tensor of the limit space.
The latter statement was announced by the second author \cite{petrunin-poly}.

Analogous statements about metric tensor and Levi-Civita connection were essentially proved by Perelman \cite{PerDC},
we only had to tie his argument with an appropriate convergence.
This part is discussed in Section~\ref{sec:DC}. 
It provides a technique that could be useful elsewhere as well.
For curvature tensor (which has a higher order of derivative), this argument cannot be extended directly;
we found a way around applying Bochner-type formulas as in \cite{petrunin-SC}.

The following statement looks like a direct corollary of the main theorem, 
and indeed, it follows from its proof but strictly speaking, it cannot be deduced directly from the main theorem alone.
We will denote by $\Sc$ the scalar curvature and $\vol^m$ the $m$-dimensional volume; that is, $m$-dimensional Hausdorff measure calibrated so that the unit $m$-dimensional cube has unit measure.

\begin{thm}{Corollary}\label{cor:Sc}
In the assumption of the main theorem,
the measures $\Sc\cdot \vol^m$ on $M_n$ weakly converge to a locally finite signed measure $\mathfrak m$  on $A$.
\end{thm}

The following subcorollary requires no new definitions.

\begin{thm}{Subcorollary}\label{cor:cor:Sc}
In the assumption of the main theorem, suppose $A$ is compact.
Then the sequence
\[s_n=\int_{M_n}\Sc \cdot\vol^m\]
converges.
\end{thm}

The main theorem in \cite{petrunin-SC} implies that if a sequence of complete $m$-dimensional Riemannian  manifolds $M_n$ has uniformly bounded diameter and uniform lower curvature bound, then 
the corresponding sequence $s_n$ is bounded;
in particular, it has a converging subsequence.
However, if $M_n$ is collapsing, then this sequence may not converge.
For example, an alternating sequence of flat 2-toruses and round 2-spheres might collapse to the one-point space; in this case, the sequence $s_n$ is $0,4\cdot\pi,0,4\cdot\pi,\dots$

From the main theorem (and the definition of weak convergence) we get the following.

\begin{thm}{Corollary}
Let $\mathfrak{K}$ be a convex closed subset of curvature tensors on $\RR^m$ 
such that 
all sectional curvatures of tensors in $\mathfrak{K}$ are at least $-1$.
Assume that $\mathfrak{K}$ is invariant with respect to the rotations of $\RR^m$.
(For example, one can take as $\mathfrak{K}$ the set of all curvature tensors with nonnegative curvature operator.)

Suppose $M_n$ is a sequence of complete $m$-dimensional Riemannian manifolds that converges to a Riemannian manifold $M$ of the same dimension.
Assume that for any $n$, all curvature tensors of $M_n$ belong to $\mathfrak{K}$, 
then the same holds for the curvature tensors of $M$.
\end{thm}

\paragraph{Remarks.}
The limit measure $\mathfrak m$ in \ref{cor:Sc} has some specific properties;
let us describe a couple of them:
\begin{itemize}
\item The measure $\mathfrak m$ vanishes on any subset of $A$ with a vanishing $(m\zz-2)$-dimensional Hausdorff measure.
In particular, $\mathfrak m$ vanishes on the set of singularities of codimension 3.
This is an easy corollary of \cite{petrunin-SC}.

\item The measure can be explicitly described on the set of singularities
of codimension~$2$.
Namely, suppose $A'\subset A$
denotes the set of all points $x$ with tangent space
$\T_xA\zz=\R^{m-2}\zz\times \Cone(\theta)$,
where $\Cone(\theta)$ is a $2$-dimensional cone
with the total angle $\theta\zz=\theta(x)\zz<2\cdot\pi$.
Then 
$$\mathfrak m|_{A'}=(2\cdot \pi-\theta)\cdot \vol^{m-2}.$$
This statement follows from \ref{prop:3parts:codim2+}.
\end{itemize}

The geometric meaning of our curvature tensor is not quite clear. 
In particular, we do not see a solution to the following problem;
compare to \cite[Conjecture~1.1]{G}.

\begin{thm}{Problem}
Suppose that the limit curvature tensor of Alexandrov space $A$ as in the main theorem has sectional curvature bounded below by $K>\kappa$.
Show that $A$ is an Alexandrov space with curvature bounded below by~$K$.
\end{thm}

The theorem makes it possible to define a curvature tensor for every \textit{smoothable} Alexandrov space.
It is expected that the same can be done for general Alexandrov space; so the following problem has to have a solution:

\begin{thm}{Problem}\label{prob:curvature}
Extend the definition of measure-valued curvature tensor to general Alexandrov spaces.
\end{thm}

If this is the case, then one may expect to have a generalization of the Gauss formula for the curvature of a convex hypersurface, which in turn might lead to a solution of the following open problems in Alexandrov geometry.
This conjecture is open even for convex sets in \textit{smoothable} Alexandrov space.

\begin{thm}{Conjecture}
The boundary of an Alexandrov space equipped with its intrinsic metric is an Alexandrov space with the same lower curvature bound.
\end{thm}

More importantly, a solution to \ref{prob:curvature} might provide nontrivial ways to deform Alexandrov space; see \cite[Section 9]{petrunin-conc}. 

\parbf{Related results.}
The result of the main theorem in dimension 2 is well known \cite[VII \S13]{AZ}.

The construction of harmonic coordinates at regular points of $\rcd$ space (in particular, Alexandrov space) given by Elia Bruè, Aaron Naber, and Daniele Semola \cite{BNS} might help to solve \ref{prob:curvature}.

The problem of introducing Ricci tensor
was studied in far more general settings \cite{G1,St,H,L}.
Curvature tensor for $\rcd$ spaces was defined by Nicola Gigli \cite{G};
it works for a more general class of spaces, but this approach does not see the curvature of singularities.
It is expected that our definitions agree on the regular locus.

\paragraph{About the proof.}
As it was stated, the 2-dimensional case is proved in \cite[VII \S13]{AZ}.
The 3-dimensional case is the main step in the proof;
the higher-dimensional case requires only minor modifications.

We subdivide the limit space $A$ into
three subsets: $A^\circ$ --- the subset of regular 
points, $A'$ --- points with singularities of codimension 2,
$A''$ --- singularities of higher codimension.
These sets are treated independently.

First, we show that limit curvature vanishes on $A''$; 
this part is an easy application of the main result in \cite{petrunin-SC}.

The $A'$-case is reduced to its partial case when the limit is isometric to the product of the real line and a two-dimensional cone.
The proof uses a Bochner-type formula (\ref{thm:bochner-formula}) and Theorem \ref{thm:extimage-of-G-and-H} which is a more exact version of the following problem from \cite{petrunin-PIGTIKAL}.

{

\begin{wrapfigure}{r}{35 mm}
\vskip-7mm
\centering
\includegraphics{mppics/pic-25}
\vskip0mm
\end{wrapfigure}

\begin{thm}{Convex-lens problem}
Let $D$ and $D'$ be two smooth discs with a common boundary that bound a convex set (a lens) $L$ in a positively-curved 3-dimensional Riemannian manifold $M$.
Assume that the discs meet at a small angle.
Show that the integral $\int_{D}k_1\cdot k_2$
is small; here $k_1$ and $k_2$ denote the principal curvatures of $D$.
\end{thm}

The $A^\circ$-case is proved by induction.
The base is the 2-dimensional case.
Further, we apply the induction hypothesis to level sets of special concave functions.
By the Gauss formula, these level sets have the same lower curvature bound. 
In the proof, we use the Bochner-type formula together with the DC-calculus developed in \cite{PerDC}.
The first step in the induction is slightly simpler.

}

As a rule, the calculus is done in the approximating sequence of Riemannian manifolds.

\parbf{Acknowledgments.}
We wish to thank Sergei Ivanov for pointing out a gap in a preliminary version of this paper,
John Lott for expressing his interest in a written version for many years,
and Alexander Lytchak for helping us to write this paper in a more readable way.
Our very special thanks to an anonymous referee, who suggested several dozens of refinements. 

The first author was partially supported by the Russian Foundation for Basic Research grant 20-01-00070, 
the second author was partially supported by the National Science Foundation grant DMS-2005279
and the Ministry of Education and Science of the Russian Federation, grant 075-15-2022-289.

\tableofcontents

\section{Formulations}

In this section, we give the necessary definitions for a precise formulation of the main theorem.
For simplicity  we will always assume that the lower
curvature bound is  $-1$;
applying rescaling, we can get the general case.

We denote by
$\Al^m$ the class of $m$-dimensional Alexandrov's spaces
with curvature $\ge -1$.

Suppose $A,A_1,A_2,\dots{} \in \Al^m$ and $A_n\GHto A$.
That is, $A_n$ converges to 
$A$ in the sense of Gromov--Hausdorff;
since $A\in \Al^m$, we have no collapse.
Denote by $a_n\:A_n\to A$ the corresponding Hausdorff approximations.
If $A$ is compact, then by Perelman's stability theorem \cite{PerStab,KapStab} we can (and will) assume that $a_n$ is a homeomorphism for every sufficiently large $n$.
In the case of noncompact limit, we assume that for any $R$, the restriction of $a_n$ to an $R$-neighborhood of the marked point is a homeomorphism to its image for every sufficiently large $n$.

We say that $A\in \Al^m$ is \emph{smoothable}
if it can be presented as a Gromov--Hausdorff limit of a non-collapsing sequence of Riemannian manifolds $M_n$ with $\sec M_n\ge-1$; here $\sec$ stands for sectional curvature.
Given a smoothable Alexandrov space $A$,
a sequence of complete Riemannian manifolds $M_n$ as above together with a sequence of approximations $a_n\:M_n\to A$
will be called \emph{smoothing} of $A$
(briefly, $M_n\smooths{} A$, or $M_n\smooths{a_n} A$).
By Perelman's stability theorem, any smoothable Alexandrov space is a topological manifold without boundary.

\subsection{Weak convergence of measures}

In this subsection, we define weak convergence of measures.
For more detailed definitions and terminology, we refer to
\cite{GMS}.

Let $X$ be a Hausdorff topological space.
Denote by $\mathfrak M(X)$ the space of signed Radon measures on $X$.
Further, denote by $C_c(X)$  the space of continuous functions on $X$
with a compact support. 

We  denote by $\langle \mathfrak m|f\rangle $ the value of $\mathfrak m\in\mathfrak M(X)$ on $f\in C_c(X)$.
We say that measures $\mathfrak m_n\in \mathfrak M(X)$ \emph{weakly converge} to $\mathfrak m\in \mathfrak M(X)$ (briefly
$\mathfrak m_n\rightharpoonup \mathfrak m$) if $\langle \mathfrak m_n|f\rangle \to \langle \mathfrak m|f\rangle $ for any $f\in C_c(X)$.

Suppose $A_n\GHto A$ with Hausdorff approximations $a_n\:A_n\zz\to A$ and
$\mathfrak m_n$ is a measure on $A_n$.
We say that $\mathfrak m_n$ \emph{weakly converges} to a measure $\mathfrak m$ on $A$ (briefly $\mathfrak m_n\rightharpoonup \mathfrak m$) if the pushforwards $\mathfrak m_n'$ of $\mathfrak m_n$ to $A$  by the Hausdorff approximations $a_n\:A_n\to A$ weakly converge to 
$\mathfrak m$.
If the condition $\langle \mathfrak m_n'|f\rangle \to \langle \mathfrak m|f\rangle $ holds only for functions $f$ with support in an open subset $\Omega\subset A$, then we say that $\mathfrak m_n$ \emph{weakly converges} to $\mathfrak m$ in $\Omega$.

Equivalently, the weak convergence can be defined using the uniform convergence of functions.
We say that  a sequence $f_n\in C_c(A_n)$
\emph{uniformly converges} to $f\in C_c(A)$
if their supports are uniformly bounded and
\[\sup_{x\in A_n}\{\,|f_n(x)-f\circ a_n(x)|\,\}\to 0.\]
Then  $\mathfrak m_n\rightharpoonup \mathfrak m$
if for any sequence $f_n\in C_c(A_n)$
with uniformly bounded supports and
uniformly converging to $f\zz\in C_c(A)$
we have $\langle \mathfrak m_n|f_n\rangle \to \langle \mathfrak m|f\rangle $.

\subsection{Test functions}\label{sec:test-functions}

In this subsection, we introduce a class of test functions and define their convergence.

Test functions form a narrow class of functions defined via a formula.
It is just one possible choice of a class containing sufficiently smooth DC functions; see the remarks in the next section.

Recall that the distance between points $x,y$ in a metric space is denoted by $|x-y|$;
we will denote by $\dist_x$ the distance function $\dist_x\:y\mapsto |x-y|$.

{\sloppy 

Suppose $A_n, A\in \Al^m$ and  $A_n\GHto A$.
Then any distance function $\dist_p\:A\zz\to\RR$ can be \emph{lifted} to $A_n$;
it means that we can choose a convergent sequence $p_n\to p$ and take the
sequence $\dist_{p_n}$.

}

Choose $r>0$ and $p\in A$.
Let us define \emph{smoothed distance function} as the average:
$$\widetilde{\dist}_{p,r} =\oint_{B(p,r)} \dist_{x}dx.$$ 
We can lift this function to
$\widetilde{\dist}_{p_n,r}\:A_n\to [0,\infty)$
by choosing some  sequence $A_n\ni p_n\to p\in A$.

We say that $f$ is a \emph{test function} if it can be expressed by the formula
$$f=\varphi( \widetilde{\dist}_{p_1,r_1}, \dots,   \widetilde{\dist}_{p_N,r_N}),$$
where $\varphi:(0,\infty)^N\to\R$ is a $C^2$-smooth function with compact support.
If for some sequences of points $A_n\ni p_{i,n}\to p_i\in A$ and $C^2$-smooth functions $\varphi_n$ that $C^2$-converge to $\varphi$ with compact support we have
$$f_n=\varphi_n( \widetilde{\dist}_{p_{1,n},r_1},\dots,   \widetilde{\dist}_{p_{N,n},r_N}),$$
then we say that $f_n$ is \emph{test-converging} to $f$ (briefly, $f_n\zz\testto f$).

\parbf{Remarks.}
Note that test functions form an algebra.

Let $M$ be a Riemannian manifold.
Note that for any open cover of $M$, there is a subordinate partition of unity of test functions.
Further, around any point of $M$ one can take a smoothed distance 
coordinate chart.
One can express any $C^2$-smooth function in these 
coordinates, and then apply partition of unity for a covering by charts.
This way, we get the following:

\begin{thm}{Claim}
On a smooth complete Riemannian manifold, test functions
include all $C^2$-smooth functions with compact support.
\end{thm}

\subsection[\texorpdfstring{$C^1$-delta convergence}{C¹-delta convergence}]%
{$\bm{C^1}$-delta convergence}\label{sec:concept}

Here we introduce $C^1$-delta convergence.
It will be necessary to formulate the main theorem in an invariant way, but, except for \ref{sec:test-convergence}, everywhere in the proofs, we will use test convergence and occasionally DC convergence instead.
(As claimed in \ref{clm:test=>smooth} test convergence implies $C^1$-delta convergence.)
By that reason, \textit{it would be wise to skip this section for the first reading}.

The $C^1$-delta convergence will be used together with other delta convergences introduced in~\ref{subsec:chart+delta}.

\parbf{Convergence of vectors.}
Let $A$ be an Alexandrov space, we denote by $\T A$ the set of all tangent vectors at all points.
So far $\T A$ is a disjoint union of all tangent cones;
let us define a convergence on it.

We will use gradient exponent $\gexp\: \T A\to A$ which is defined in \cite{AKP}.
Given a vector $V\in \T A$, it defines its radial curve $\gamma_{V}\:t\zz\mapsto \gexp (t\cdot V)$.
We say that a sequence of vectors $V_n\in \T A$ \emph{converges} to $V\in \T A$ (briefly, $V_n\to V$) if $\gamma_{V_n}$ converges to $\gamma_V$ pointwise.
Since the radial curve $\gamma_V$ is $|V|$-Lipschitz, we get that any bounded sequence of vectors with base points in a bounded set has a converging subsequence of $\gamma_{V_n}$.
Further, the pointwise limit of such curves is a radial curve as well.
Therefore, any bounded sequence of tangent vectors with base points in a bounded set has a converging sequence.

In a similar fashion, we can define the convergence of tangent vectors to sequences of Alexandrov spaces $A_n$ that converge to $A$.
That is, if $V_n\in \T A_n$ is a bounded sequence of tangent vectors at points on a bounded distance to the base points, then it has a subsequence that converges to some vector $V\in \T A$.

Note that 
\[|V|\le \liminf_{n\to\infty} |V_n|\]
and the inequality might be strict.

Recall that if $V\in \T_p$ is the unit vector in the direction of $[pq]$, then $\gamma_V$ is a unit-speed parametrization of $[pq]$.
Using this we get the following observation;
it provides a way to apply the convergence.

\begin{thm}{Observation}\label{obs:unique-geod}
Let $M_n\smooths{} A$ be a smoothing, $p_n,q_n\in M_n$, and $p_n\to p$, $q_n\to q$ as $n\to \infty$.
Denote by $V_n\in \T_{p_n}$ and $V\in \T_p$ the directions of geodesics $[p_nq_n]$ and $[pq]$.
Suppose that there is a unique geodesic $[pq]$ in $A$.
Then $V_n\zz\to V$.
\end{thm}

\parbf{$\bm{C^1}$-delta smoothness.}
Given a function $f\:A\to \RR$ and a vector $V\in \T A$, set
\[Vf=(f\circ\gamma_V(t))'|_{t=0}.\]
Note that $Vf$ is defined for all DC functions and, in particular, all test functions.

Two vectors $V,W\in \T_pA$ will be called \emph{$\delta$-opposite} if
$1-\delta< |V|\le 1$,
$1-\delta< |W|\le 1$,
and $|\langle X,V\rangle +\langle X,W\rangle|<\delta$ for any unit vector $X\in \T_p A$.
We say that $V,W\in \T_pA$ are opposite if they are $\delta$-opposite for any $\delta>0$;
in this case, they are both unit vectors and make angle $\pi$ to each other.

A function $f\:A\to\RR$ is called $C^1$-delta smooth if for any compact set $K\subset A$ and $\eps>0$ there is $\delta>0$ such that any sequence of points $p_n\to p\in K$ and unit vectors $V_n\in \T_{p_n} A$ that converges to a vector $V\in \T_p A$ that has a $\delta$-opposite vector we have
\[|Vf-\acute{\lim_{n\to\infty}} V_nf|<\eps,\]
where ``\,$\acute{\lim}$\,'' stands for an arbitrary partial limit.

Suppose $M_n\smooths{} A$.
A sequence of $C^1$-smooth functions $f_n\:M_n\zz\to\RR$ is called \emph{$C^1$-delta converging} to $f\:A\to \RR$ (briefly, $f_n\smoothto f$) if $f_n$ converges to $f$ pointwise and
for any compact set $K\subset A$ and any $\eps>0$ there is $\delta>0$
such that if a sequence of unit vectors $V_n\in \T_{p_n} M_n$ converges to a vector $V\zz\in \T_p A$ such that $p\in K$ and $V$ has a $\delta$-opposite vector, then we have
\[|Vf-\acute{\lim_{n\to\infty}} V_nf_n|<\eps.\]

\begin{thm}{Claim}\label{clm:test=>smooth}
Any test function is $C^1$-delta smooth.
Moreover, for any smoothing $M_n\smooths{} A$,
sequence of test functions $f_n\:M_n\zz\to\RR$,
and a test function $f\:A\zz\to\RR$ we have
\[f_n\testto f
\qquad\Longrightarrow\qquad
f_n\smoothto f.\]

\end{thm}

\parit{Proof.}
Let $V$ and $W$ be $\delta$-opposite vectors in $\T_pA$.
Note that for almost all points $q\in A$, we have 
\[|V\dist_{q}+W\dist_{q}|<\delta.\]
It follows that 
\[|V\widetilde{\dist}_{q,r}+W\widetilde{\dist}_{q,r}|<\delta\eqlbl{eq:|(V+W)tildedist|<delata}\]
for any $q\in A$ and $r>0$.

Suppose $V_n$ is a sequence of unit tangent vectors on $M_n$ such that $V_n\to V$; that is, $\gamma_{V_n}\to \gamma_V$ as $n\to\infty$.
By monotonicity of radial curves \cite[16.32]{AKP}, we get 
\[V\dist_{q}\le  \liminf_{n\to\infty} V_n\dist_{q_n}\]
if $q_n\to q$.
Integrating, we get 
\[V\widetilde{\dist}_{q,r}\le  \liminf_{n\to\infty}V_n\widetilde{\dist}_{q_n,r}.\]

Suppose $V$ has a $\delta$-opposite vector $W$.
We can assume that $W$ is a unit geodesic vector; that is, there is a geodesic $[ps]$ in the direction of $W$.
Moreover, we can assume that $[ps]$ is a unique geodesic from $p$ to $s$.
Choose points $s_n$ and $p_n$ that converge to $s$ and $p$ respectively.
By \ref{obs:unique-geod}, the directions $W_n$ of $[p_ns_n]$ converge to $W$.
Note that $W_n$ is $\delta$-opposite to $V_n$ for all large $n$.

Repeating the above argument, we get 
\[W\widetilde{\dist}_{q,r}\le  \liminf_{n\to\infty}W_n\widetilde{\dist}_{q_n,r}.\]
Applying \ref{eq:|(V+W)tildedist|<delata} we get $C^1$-delta convergence of $\widetilde{\dist}_{q_n,r}$ and, in particular, $C^1$-delta smoothness of $\widetilde{\dist}_{q,r}$.
Applying the definition of test function, we get the result.
\qeds

Recall (\ref{sec:test-functions}) that for any smoothing $M_n\smooths{} A$ and a test function $f\:A\to\RR$ there are test functions  $f_n\:M_n\to \RR$ such that $f_n\testto f$.

\begin{thm}{Corollary}\label{cor:liftability}
Given a smoothing $M_n\smooths{} A$ and a test function $f\:A\to\RR$,
there is a sequence of $C^1$-smooth functions $f_n\:M_n\to \RR$ such that $f_n\smoothto f$.
\end{thm}

\parbf{Remarks.}
In the next section, we define measure-valued tensor as a functional on an array of test functions.
Note that, one test function might have very different presentations that lead to different test convergences.
Thus to prove the invariance of measure-valued curvature tensor we need to use the $C^1$-delta convergence which is more general than test convergence.
We could use other classes of functions as well. 
For example, a subclass of DC$_0$ functions (see Section \ref{sec:DC})
or a subclass of $C^1$-delta function (see Section \ref{sec:concept}).
Of course, we have to have an analog of \ref{cor:liftability} for the chosen class.
Hope a more natural setting will be found eventually.

\subsection{Tensors}\label{subsec:tensors}

In this subsection, we define measure-valued tensors on Alexandrov spaces.
Basically, we reuse the derivation approach to vector fields in classical differential geometry.
This definition will be used in Claim \ref{clm:weak-partial-limit} that reduces the main theorem to Proposition~\ref{prop:3parts} and will not show up ever after.

Let $A\in \Al^m$.
Recall that $\mathfrak M(A)$
denotes the space of signed Radon measures on $A$.
A \emph{measure-valued vector field} $\mathfrak{v}$ on $A$
is a linear map that takes a test function,
spits a measure in $\mathfrak M(A)$,
and satisfies the \emph{chain rule}:
for any collection of test functions $f_1,\dots,f_k$
and a $C^2$-smooth function $\phi\:\RR^k\to\RR$, we have
$$\mathfrak{v}(\phi(f_1,\dots,f_n))
=
\sum_{i=1}^n (\partial_i\phi)(f_1,\dots,f_n)\cdot\mathfrak{v}(f_i).
$$

In the same way, we define (contravariant) measure-valued tensor fields.
Namely, a \emph{measure-valued tensor field} $\mathfrak{t}$ of valence $k$ on $A$ is a multilinear map that takes a $k$-array of test  functions, spits a measure in $\mathfrak M(A)$, and satisfies the chain rule in each of its arguments.

Suppose that $x_1,\dots,x_m$ are local coordinates in an $m$-dimensional Riemannian manifold~$M$.
Then a measure-valued vector field $\mathfrak{v}$ on $M$ can be described by $m$ components, these are measures $(\vv(x_1),\zz\dots,\vv(x_m))$;
these components transform by contravariant rule under change of coordinates.

By the definition of measure-valued vector field, we get
\[\mathfrak{v}(f)=\sum_{i}\partial_i f\cdot \vv(x_i).\]
Similarly, for arbitrary $k$, a measure-valued tensor field of valence $k$ is defined by $m^k$ components 
$\mathfrak{t}(x_{i_1},\dots,x_{i_k})$; namely,
\[\mathfrak{t}(f_1,\dots,f_k)
=
\sum_{i_1,\dots,i_k}
\partial_{i_1} f_1 
\cdots 
\partial_{i_k} f_k
\cdot \mathfrak{t}(x_{i_1},\dots,x_{i_k}).\]

Note that if $T$ is a smooth contravariant tensor field then $\mathfrak{t}=T\cdot \vol$ is a measure-valued tensor field.
In other words, usual tensor fields might be considered as a subspace of measure-valued tensor fields.

\begin{rdef}{Definition}
Let $M_n\smooths{} A$ be a smoothing.
Assume that $\mathfrak{t}_n$ is a sequence of 
 measure-valued tensor fields on $M_n$  and $\mathfrak{t}$ is a measure-valued tensor field on $A$,
all of the same valence~$k$.
We say that $\mathfrak{t}_n$ \emph{weakly converges} to  $\mathfrak{t}$
(briefly $\mathfrak{t}_n\rightharpoonup\mathfrak{t}$) if
\[f_{i,n}\smoothto f_i\quad\text{for all\ } i
\quad\Longrightarrow\quad
\mathfrak{t}_n(f_{1,n},\dots,f_{k,n})\zz\rightharpoonup\mathfrak{t}(f_{1},\dots,f_{k})\]
for arbitrary $k$ sequences $f_{1,n},\dots,f_{k,n}$ of $C^1$-smooth functions and test functions $f_1,\dots,f_k\:A\zz\to\RR$.
\end{rdef}

\subsection{Dual curvature tensor}

The curvature of Riemannian manifold $M$ is usually described by a tensor of valence 4 that will be denoted by $\Rm$.
We will use a \emph{dual curvature tensor} --- 
a curvature tensor written in a dual form that will be denoted by $\Qm$;
it is a tensor field of valence $2\cdot(m-2)$ defined the following way:
\begin{align*}
\Qm(X_1,&\dots,X_{m-2},Y_1,\dots,Y_{m-2})
= 
\\
&=
\Rm(*(X_1\wedge\dots\wedge X_{m-2}),*(Y_1\wedge\dots\wedge Y_{m-2})),
\end{align*}
where $X_i,Y_i$ are vector fields on $M$ and  ${*}\:(\bigwedge^{m-2}\T)M\zz\to(\bigwedge^2\T)M$ is the  Hodge star operator.
This definition will be used further mostly for gradient vector fields of semiconcave functions.

In addition, we will need a measure-valued version of $\Qm$ denoted by $\qm$;
it will be called \emph{dual measure-valued curvature tensor}.
Namely, we define 
\[\qm(f_1,\zz\dots,f_{m-2},g_1,\zz\dots,g_{m-2})\]
as the measure with density
\[\Qm(\nabla f_1,\dots,\nabla f_{m-2},\nabla g_1,\dots,\nabla g_{m-2}): M\to\R.\]

\parbf{Remarks.}
Note that 
$$\Qm(X_1,\dots,X_{m-2},X_1,\dots,X_{m-2})
=
|X_1\wedge\dots \wedge X_{m-2}|^2\cdot K_\sigma, $$
where $K_\sigma$ is the sectional curvature of $M$ 
on a plane $\sigma$ orthogonal to $(m-2)$-vector
$X_1\wedge\zz\dots \wedge X_{m-2}$.
Hence, the sectional curvatures of $M$ and therefore its curvature tensor $\Rm$ can  be computed from
$\qm$.
By the symmetry
\begin{align*}
\qm(f_1,&\dots,f_{m-2},g_1,\dots,g_{m-2})
\\
&=
\qm(g_1,\dots,g_{m-2},f_1,\dots,f_{m-2}),
\end{align*}
the density of $\qm$ is defined by the sectional curvature.
Therefore measure-valued tensor $\qm$
gives an equivalent description of the curvature of Riemannian manifolds.

As you will see further, the described dual form of curvature tensor behaves better in the limit;
in particular, it makes it possible to formulate \ref{prop:3parts:codim2+}.

In the 2-dimensional case, the valence of $\qm$ is $0$;
in this case, $\qm$ coincides with the curvature measure --- the standard way to describe the curvature of surfaces, 
\cite{Resh,AZ}.
For a smooth surface, the density of this curvature measure with respect to the area
is its Gauss curvature.
In this case, it is known that \textit{curvature measures are stable under smoothing} \cite[VII \S13]{AZ};
in other words, our main theorem is known in the two-dimensional case.

\subsection{Formulation and plan}
  
\begin{thm}{Main theorem}\label{main}
Consider a smoothing $M_n\smooths{} A$.
Denote by $\qm_n$ the dual measure-valued curvature tensor on $M_n$.
Then there is a measure-valued tensor $\qm$ on $A$ such that $\qm_n\rightharpoonup \qm$.
\end{thm}

Let $A$ be an $m$-dimensional Alexandrov space without boundary. 
Let us partition $A$ into three subsets $A^\circ$, $A'$, and $A''$:
\begin{itemize}
\item $A^\circ$ is the set of regular points in $A$; that is, the set of points with tangent cone isometric to the Euclidean space.
\item $A'$ --- the set of points in $A\backslash A^\circ$ with an isometric copy of $\RR^{m-2}$ in their tangent space;
in other words, for any $p\in A'$, the tangent space $\T_p$ is isometric to the product $\Cone(\theta)\zz\times\RR^{m-2}$ where $\Cone(\theta)$ denotes a two-dimensional cone with the total angle $\theta\zz=\theta(p)<2\cdot \pi$.
\item $A''$ --- the remaining set; this is the set of points with tangent space that does \textit{not} contain an isometric copy of $\RR^{m-2}$.
\end{itemize}
According to \cite{li-naber}, $A'$ is countably $(m-2)$-rectifiable,
and $A''$ is countably $(m-3)$-rectifiable. 

Observe that the set of regular points $A^\circ$ can be presented as
$$A^{\circ}=\bigcap_{\delta>0} A^\delta,$$
where $A^\delta$ denotes the set of $\delta$-strained points of $A$.

Let $M$ be an $m$-dimensional Riemannian manifold.
Denote by $K_{max}(x)$ the maximal sectional curvature at $x\in M$.
The following statement is a direct corollary of the main result in \cite{petrunin-SC}:

\begin{thm}{Corollary}\label{cor:Kmax}
Given an integer $m\ge 0$, there is a constant $\const(m)$ such that the following holds:

Let $M$ be an $m$-dimensional Riemannian manifold 
(possibly noncomplete)
with sectional curvature bounded below by $-1$.
If for some $r<1$ the closed ball $\bar B(p,2\cdot r)_M$ is compact,
then 
$$\int_{B(p,r)_M} K_{max}\le \const(m)\cdot r^{m-2}.$$

\end{thm}

\begin{thm}{Observation}\label{q=<K}
There is another constant $\const'(m)$ such that 
\begin{align*}
|\Qm(X_1,&\dots,X_{m-2},Y_1,\dots,Y_{m-2})|\le
\\
&\le 
\const'(m)\cdot K_{max}\cdot|X_1\wedge\dots\wedge X_{m-2}|\cdot |Y_1\wedge\dots\wedge Y_{m-2}|
\end{align*}

\end{thm}

Note that \ref{cor:Kmax} and \ref{q=<K} imply the following:

\begin{thm}{Claim}\label{clm:weak-partial-limit}
Given a smoothing $M_n\smooths{} A$,
test functions $f_i\:A\zz\to \RR$,
and sequences of $C^1$-smooth functions $f_{i,n}\:M_n\to \RR$ such that 
$f_{i,n}\zz\smoothto f_i$ the sequence of measures 
$\qm_n(f_{1,n},\zz\dots f_{2\cdot m-4,n})$ has a weakly converging subsequence.

Moreover, the subsequence can be chosen simultaneously for several choices of function arrays so that it meets the chain rule.
More precisely, choose $i$; fix all functions $f_{1,n},\zz\dots, f_{2\cdot m-4,n}$ except $f_{i,n}$;
suppose $\hat\qm_{i,n}(f_{i,n})=\qm_n(f_{1,n},\zz\dots, f_{2\cdot m-4,n})$.
Assume $h_{j,n}\:M_n\to \RR$ are $C^1$-smooth functions such that 
$h_{j,n}\zz\smoothto h_j$, each $h_j$ is a test function, and 
\[h_{0,n}=\phi(h_{1,n},\dots,h_{k,n})\]
for a fixed $C^2$-function $\phi\:\RR^k\to\RR$.
Then the sequence of measure arrays $\hat\qm_{i,n}(h_{0,n}),\zz\dots,\hat\qm_{i,n}(h_{k,n})$
has a partial limit $\hat\qm_{i}(h_0),\zz\dots,\hat\qm_{i}(h_k)$, and
\[\hat\qm_{i}(h_0)=\sum_{j=1}^k (\partial_j\phi)(h_1,\dots,h_k)\cdot\hat\qm_{i}(h_j).\]

\end{thm}

By Perelman's stability theorem, the space $A$ in the claim is a topological manifold.
In particular, $A$ has no boundary;
in other words, the singular set in $A$ has codimension at least~2.
Together with the claim, it implies that the main theorem (\ref{main}) follows from the next statement.

\begin{thm}{Proposition}\label{prop:3parts}
Let $M_n\smooths{} A$ and $\dim A=m$.
Suppose $h_{1},\zz\dots, h_{m-2}$ are 
test functions on $A$ and
$h_{1,n},\dots, h_{m-2,n}$ are 
$C^1$-smooth functions on $M_n$ 
such that $h_{i,n}\zz\smoothto  h_i$ for each $i$.
Let $\mathfrak m_1$ and $\mathfrak m_2$ be two measures on $A$ that are weak partial limits of the sequence of measures $\qm_n( h_{1,n},\zz\dots,  h_{m,n}, h_{1,n},\zz\dots,  h_{m,n})$ on $M_n$.
Then the following statements hold:
\begin{enumerate}[(i)]
\item\label{prop:3parts:codim3} $\mathfrak m_1|_{A''}=\mathfrak m_2|_{A''}=0$;

\item\label{prop:3parts:codim2} $\mathfrak m_1|_{A'}=\mathfrak m_2|_{A'}$;

\item\label{prop:3parts:reg} $\mathfrak m_1|_{A^\circ}=\mathfrak m_2|_{A^\circ}$.
\end{enumerate}
\end{thm}

The three parts of the proposition will be proved below in sections \ref{sec:codmi=3}, \ref{sec:codmi=2}, and \ref{sec:ref} respectively.

\part*{Proofs}
\addcontentsline{toc}{part}{Proofs}

\section{Singularities of codimension 3}\label{sec:codmi=3}

\parit{Proof of \ref{prop:3parts}(\ref{prop:3parts:codim3}).}
According to \cite[10.6]{BGP}, $A''$ has a vanishing $(m-2)$-dimensional Hausdorff measure;
that is, $A''$ can be covered by a countable family of balls $B(x_i,r_i)$ such that $\sum r_i^{m-2}$ is arbitrarily small.
Therefore, \ref{q=<K} and \ref{cor:Kmax} imply the statement.
\qeds

\section{Singularities of codimension 2}\label{sec:codmi=2}

The following lemma will be proved in Section~\ref{sec:bilip}.

\begin{thm}{Lemma}\label{lem:A-prime-Q}
Let $A$ be an $m$-dimensional Alexandrov space without boundary.
Then the subset $A'\subset A$ can be covered by a countable set of compact sets $Q_i$ such that each $Q_i$ admits a bi-Lipschitz embedding into $\RR^{m-2}$.
\end{thm}

Let $\bm{h}\:A\to \RR^k$ be a Lipschitz map defined on an $m$-dimensional Alexandrov space without boundary.
Suppose $Q\subset A$ is a closed subset such that there is a bi-Lipschitz embedding $s\:Q\to\RR^k$.
By the generalized Rademacher theorem, the metric differential of $s^{-1}$ is defined almost everywhere in the domain of definition of $s^{-1}$.
Moreover, the metric differential is defined by a bilinear form; 
its determinant is the Jacobian of $s^{-1}$, briefly $\jac s^{-1}$.
The same way we can define $\jac (\bm{h}\circ s^{-1})$ (we can apply the standard Rademacher theorem this time).
Further, set $\jac(\bm{h}|_Q)\zz=\jac (\bm{h}\circ s^{-1})/\jac s^{-1}$.
It is straightforward to check that this definition is $\vol^k$-almost-everywhere independent of the choice of $s$.

Consider the function
\[\theta(p)=2\cdot\pi\cdot \tfrac{\vol^{m-1} \Sigma_p}{\vol^{m-1} \mathbb{S}^{m-1}},
\eqlbl{eq:theta}\]
where $\Sigma_p$ denotes the space of directions at $p$.
According to \cite[7.14]{BGP}, $\theta\:A\to \RR$ is lower-semicontinuous.

Note that $\theta$ is identically $2\cdot\pi$ on $A^\circ$.
Further note that for any point $p\in A'$, its tangent cone is isometric to the product space 
$\Cone(\theta)\times\RR^{m-2}$, where $\theta=\theta(p)<2\cdot\pi$. 
Since $\vol^{m-2}(A'')=0$, the measure $(2\pi-\theta)\cdot\vol^{m-2}$ vanishes on $A''$.

Note that \ref{prop:3parts}\textit{(\ref{prop:3parts:codim2})} follows from Lemma~\ref{lem:A-prime-Q} and the following statement;
it will be proved in \ref{subsec:3d}--\ref{subsec:4d}.

\begin{thm}{Proposition}\label{prop:3parts:codim2+}
Let $\mathfrak m$ be one of two limit measures $\mathfrak m_i$ in \ref{prop:3parts} and $\bm{h}\zz=(h_1,\zz\dots,h_{m-2})\:A\zz\to \RR^{m-2}$ is an array of test functions.
Suppose that $Q\subset A$ is a compact subset that admits a bi-Lipschitz embedding into $\RR^{m-2}$.
Then
\[\mathfrak m|_{Q}=(2\cdot\pi-\theta)\cdot (\jac(\bm{h}|_Q))^2\cdot \vol^{m-2}.\]
\end{thm}

\subsection{Gauss and mean curvature estimates}

\begin{thm}{Theorem}\label{thm:extimage-of-G-and-H}
Let $f$, $h$ be a pair of strongly convex smooth 1-Lipschitz functions defined on an open set of a 3-dimensional Riemannian manifold.
Suppose that
\begin{enumerate}[(i)]
\item $|\nabla f|\ge 1$ and
\[|\nabla (f+h)|<\eps\cdot|\nabla f|\] 
for some fixed positive $\eps<\tfrac12$;
\item for some $a,b\in \RR$, the set
\[W_{a,b}=\set{p\in M}{f(p)=a,\  h(p)\le b}\]
is compact.
\end{enumerate}
Denote by $k_1(p)\le k_2(p)$, 
$H(p)\zz=k_1(p)+ k_2(p)$
and
$G(p)\zz=k_1(p)\cdot k_2(p)$, the principal, mean, and Gauss curvatures of $W_{a,b}$ at $p$.
Then
\[\int_{W_{a,b}}G\le 100\cdot\eps
\eqlbl{intG<100eps}\]
and 
\[\int_{W_{a,b}}H\le 10\cdot \sqrt{\eps}\cdot \length(\partial{W_{a,b}}).
\eqlbl{intH<eps}\]
\end{thm}

The proof is based on the 2-dimensional case of the following statement,
which is the integral Bochner formula with Dirichlet boundary condition.

\begin{thm}{Proposition}\label{prop:bochner-dirichle}
Assume $\Omega$ is a compact domain with smooth boundary $\partial \Omega$ in a Riemannian manifold
and $f$ is a smooth function that vanishes on $\partial \Omega$.
Then
\[\int\limits_\Omega \left(|\Delta f|^2
-|\Hess f|^2
-\langle\mathrm{Ric}(\nabla f),\nabla f\rangle\right)
=\int\limits_{\partial\Omega}
H\cdot|\nabla f|^2,\]
where $H$ denotes the mean curvature of $\partial \Omega$.
\end{thm}

\parit{Proof of \ref{thm:extimage-of-G-and-H}.}
Equip $W_{a,b}$ with unit normal vector field $n\zz=\tfrac{\nabla f}{|\nabla f|}$.
Let 
\[S_p\:\T_pW_{a,b}\zz\to \T_pW_{a,b}\]
be the corresponding shape operator, so $S_p\:v\mapsto\nabla_vn$.
Since $f$ is strongly convex, we have that 
\[\langle S_p(v),v\rangle\ge \delta\cdot|v|^2\]
for a fixed value $\delta>0$ and any tangent vector $v\in \T_pW_{a,b}$. 

Note that the restriction $u=h|_{W_{a,b}}$ is strongly convex.
Moreover, 
\[\Hess_pu(v,v)\ge (1-\eps)\cdot \langle S_p(v),v\rangle\eqlbl{Hess=<shape}\]
for any $p\in W_{a,b}$ and $v\in\T_pW_{a,b}$.
Indeed, consider the geodesic $\gamma$ in $W_{a,b}$ such $\gamma(0)=p$ and $\gamma'(0)=v$.
Set $w\zz=\gamma''(t)$.
Note that 
\begin{align*}
w&=-\langle S_p(v),v\rangle \cdot n,
\intertext{Since $h$ is strongly convex, $\Hess_p h\ge 0$; therefore}
(\Hess_pu)(v,v)&=(\Hess_p h)(v,v)+\langle \nabla_p h,w\rangle\ge
\\
&\ge-\tfrac{\langle\nabla_p h,\nabla_p f\rangle}{|\nabla_p f|}\cdot\langle S_p(v),v\rangle\ge
\\
&\ge (1-\eps)\cdot|\nabla_p f|\cdot\langle S_p(v),v\rangle.
\end{align*}
Since $|\nabla f|\ge 1$, \ref{Hess=<shape} follows.

Since $\langle S_p(v),v\rangle\ge 0$ and $\eps<\tfrac12$, the inequality \ref{Hess=<shape} implies that 
\[4\cdot \det(\Hess_pu)\ge G(p)
\eqlbl{eq:det>=G}\]
and
\[-2\cdot \Delta u\ge H(p)
\eqlbl{eq:trace>=H}\]
for any $p\in W_{a,b}$.

Denote by  $\lambda_1(p),\lambda_2(p)$ the eigenvalues of  $\Hess_p u$, so
\begin{align*}
\trace(\Hess u)&=\Delta u=\lambda_1+\lambda_2,
\\
|\Hess u|^2&=\lambda_1^2+\lambda_2^2,
\\
\det(\Hess u)&=\lambda_1\cdot\lambda_2,
\intertext{and hence}
2\cdot\det(\Hess u)
&=|\Delta u|^2
-|\Hess u|^2.
\end{align*}

Since $W_{a,b}$ is two-dimensional, by Proposition~\ref{prop:bochner-dirichle} we get that
\[\int\limits_{W_{a,b}} 
2\cdot\det(\Hess u)
=\int\limits_{W_{a,b}} 
K\cdot|\nabla u|^2
+
\int\limits_{\partial W_{a,b}}
\kappa\cdot|\nabla u|^2,\]
where $\kappa\ge 0$ is the geodesic curvature of $\partial W_{a,b}$
and $K$ is the curvature of $W_{a,b}$.

Since $u$ is a convex function that vanishes on the boundary of $W_{a,b}$,
it has a unique critical point, which is its minimum.
By the Morse lemma,  $W_{a,b}$ is a disc.
Therefore, by the Gauss--Bonnet formula, we get that
\[\int_{W_{a,b}} K+\int_{\partial{W_{a,b}}}\kappa=2\cdot\pi.\]
Whence,
\[\int\limits_{W_{a,b}} 
\det(\Hess u)
\le\pi\cdot\sup_{p\in{W_{a,b}}}|\nabla_p u|^2.\]

Note that $\nabla_p u$ is the projection of $\nabla_ph$ to $\T_pW_{a,b}$.
Therefore,
\begin{align*}
|\nabla_p u|^2&=|\nabla_p h|^2-\langle\nabla_p h,n\rangle^2\le
\\
&\le1-(1-\eps)^2<
\\
&<2\cdot\eps.
\end{align*}
It follows that 
\[\int\limits_{W_{a,b}} 
\det(\Hess u)
\le2\cdot \pi\cdot\eps.\]
Applying \ref{eq:det>=G}, we obtain \ref{intG<100eps}.

Similarly,  by the divergence theorem, we get that
\[-\int\limits_{W_{a,b}} \Delta u=\int\limits_{\partial{W_{a,b}}} |\nabla u|.\]
Whence \ref{eq:trace>=H} implies 
\[\int\limits_{W_{a,b}} H\le 10\cdot \sqrt{\eps}\cdot \length(\partial{W_{a,b}}).\]
\qeds

\subsection{Curvature of level sets}

Let $M$ be a 3-dimensional Riemannian manifold.
Choose a smooth function $f\:M\to\RR$.
Consider its level sets 
\[L_c=\set{x\in M}{f(x)=c}.\]
If the level set $L_c$ is a smooth surface in a neighborhood of $x\in L_c$,
then denote by $k_1(x)\le k_2(x)$ the principal curvatures of $L_c$ at $x$.
In this case, set 
\begin{align*}
G(x)&=k_1(x)\cdot k_2(x),
\\
H(x)&=k_1(x)+ k_2(x);
\end{align*}
that is, $G(x)$ and $H(x)$ are Gauss and mean curvature of $L_c$ at $x$.

Recall that $\Cone(\theta)$ denotes a 2-dimensional cone with the total angle $\theta$.

\begin{thm}{Theorem}\label{thm:HG-converge}
Let $M_n\smooths{}\Cone(\theta)\times \RR$ and $f_n\: M_n\to \RR$ be a sequence of strongly concave smooth 1-Lipschitz functions.
Suppose that $\sec M_n\ge -\tfrac1n$ for each $n$, and $f_n$ converges as $n\to \infty$ to the $\RR$-coordinate $f\:(x,t)\mapsto t$ on $\Cone(\theta)\times \RR$.
Then $G_n$ and $H_n$ (the Gauss and mean curvatures of the level sets of $f_n$) weakly converge to zero.
\end{thm}

\parit{Proof.}
Choose $p\in \Cone(\theta)\times \RR$; set $a=f(p)$.

By the theorem of Artem Nepechiy \cite{Nepechiy},
there is a $(-2)$-concave function $\rho$ defined in an $r$-neighborhood of $p$ such that $\rho(x)=-|p-x|^2+o(|p-x|^2)$.
Moreover, the function $\rho$ is \emph{liftable};
that is, there is a sequence of $(-2)$-concave $\rho_n\:M_n\to\RR$ that converges to $\rho$.

Consider a point $q\in \Cone(\theta)\times \RR$ \textit{above} $p$; that is, its $\RR$-coordinate is larger, and its $ \Cone(\theta)$-coordinate is the same.
If the $\RR$-coordinate of $q$ is large, then $\dist_q+f$ is $\lambda$-concave for small $\lambda>0$ and it has a nonstrict minimum at $p$.
Therefore, given $\lambda>0$, we can find $q$ so that the sum $s=f+\dist_q+\lambda\cdot \rho$ is $(-\lambda)$-concave and has a strict maximum at $p$.
Moreover
\[-\tfrac12\cdot\lambda\cdot|p-x|_{\Cone(\theta)\times \RR}^2\ge s(x)-s(p)\ge -\tfrac32\cdot\lambda\cdot|p-x|_{\Cone(\theta)\times \RR}^2\]
and therefore
\[|\nabla_xs|\le 10\cdot\lambda\cdot|p-x|_{\Cone(\theta)\times \RR}
\eqlbl{O(p-x)}\]
if $|p-x|_{\Cone(\theta)\times \RR}$ is sufficiently small; say if $|p-x|_{\Cone(\theta)\times \RR}\zz\le \tfrac r{10}$.

Choose a sequence of points $q_n\in M_n$ that converges to $q$.
Let us apply the Green--Wu smoothing procedure to $\dist_{q_n}+\lambda\cdot \rho_n$;
denote by $h_n$ the obtained function; we can assume that $|h_n-\dist_{q_n}-\lambda\cdot \rho_n|\to 0$.
Observe that \ref{O(p-x)} implies that the first condition in \ref{thm:extimage-of-G-and-H} is met for all large $n$ in an $\tfrac r2$-neighborhood of $p_n$ with $\eps=10\cdot\lambda\cdot r$.
Moreover, one can choose $b$ so that the second condition is satisfied and $B_n=B(p_n,r/10)\cap L_a\zz\subset W_{a,b}$.
Applying \ref{thm:extimage-of-G-and-H}, we get that for any $\delta>0$, we have 
\[
\int_{B_n}G_n<\delta,
\quad\text{and}\quad
\int_{B_n}H_n<\delta.
\]
for all large $n$.
It remains to integrate the obtained inequalities by $a$ and pass to a limit as $n\to\infty$.
\qeds

For a product $\RR^{m-2}\times \Cone(\theta)$, denote by $\mathcal{V}$ the $\vol^{m-2}$-measure on the vertical line $\RR^{m-2}\zz\times \{0\}$.
Further, consider the \emph{curvature measure} 
\[\omega=(2\cdot\pi-\theta)\cdot \mathcal{V}\]
on $\RR^{m-2}\times \Cone(\theta)$.

\begin{thm}{Corollary}\label{cor:Ricci}
Suppose that $M_n\smooths{}  \Cone(\theta)\times \RR$ and $f_n\: M_n\to \RR$ be as in \ref{thm:HG-converge}.
Set $u_n\zz=\tfrac{\nabla f_n}{|\nabla f_n|}$.
Then 

\begin{enumerate}[(i)]
\item\label{cor:Ricci:Ricci} $\langle \Ric u_n,u_n\rangle$ weakly converges to zero.
\item\label{cor:Ricci:vw}
Let $v_n$ and $w_n$ be sequences of uniformly bounded, continuous vector fields on $M_n$.
Suppose that $\langle v_n,u_n\rangle$ and $\langle w_n,u_n\rangle$ converge uniformly as $n\to \infty$ to some constants $a$ and $b$ respectively.
Then 
\[\Qm(v_n,w_n)\rightharpoonup a\cdot b\cdot \omega,\]
where $\omega$ is the measure on $\Cone(\theta)\times \RR$  described above.

\end{enumerate}

\end{thm}

\parit{Proof; (\ref{cor:Ricci:Ricci}).}
Passing to a subsequence if necessary, we can assume weak convergence of $\langle \Ric u_n,u_n\rangle\zz\cdot\vol^3$ to a measure $\mathfrak m$ on $\Cone(\theta)\times \RR$.
Since $\sec M_n\ge -\tfrac1n$, we have that $\mathfrak m\ge 0$.
Therefore it is sufficient to show that $\mathfrak m\le 0$.

By \ref{thm:bochner-formula} we have that the following equality
\begin{align*}
\int\limits_\Omega &\phi_n\cdot \langle \Ric u_n,u_n\rangle =
\\
&=
\int\limits_\Omega [\phi_n \cdot G_n+H_n\cdot\<u_n,\nabla\phi_n\>-\<\nabla\phi_n,\nabla_{u_n} {u_n}\>]
\end{align*}
holds for any function $\phi_n$ with compact support on $M_n$,
assuming that all expressions in the formula have sense.

It remains to find a sequence of nonnegative functions $\phi_n\:M_n\to \RR$ with compact support that converges to a $\phi\:\Cone(\theta)\times \RR\to \RR$ such that (1) $\phi$ is unit in a neighborhood of a given point $p\in \Cone(\theta)\times \RR$ and (2) we have control on the three terms on the right-hand side of the formula; the latter means that we have the following weak convergences:
\[
\begin{aligned}
\phi_n \cdot G_n&\rightharpoonup 0,
&
H_n\cdot\<u_n,\nabla\phi_n\>&\rightharpoonup 0,
&
\<\nabla\phi_n,\nabla_{u_n} {u_n}\>&\rightharpoonup 0.
\end{aligned}
\eqlbl{eq:3->0}
\]

For the first convergence, it is sufficient to choose the sequence $\phi_n$ so that in addition it is universally bounded.
Indeed, since $|\nabla f_n|\to 1$, we have that \ref{thm:HG-converge}
implies the first convergence in~\ref{eq:3->0}.

Similarly, to prove the second convergence in \ref{eq:3->0}, it is sufficient to assume in addition that $|\nabla\phi_n|$ is universally bounded and apply \ref{thm:HG-converge}.

To prove the last convergence in \ref{eq:3->0}, note that
$|\nabla_{u_n} u_n|\zz\rightharpoonup 0$ away from the singular locus.
The latter follows from \ref{lem:test-delta} and \ref{lem:test-delta-partial-g}.
Indeed, $\nabla_{u_n} u_n$ can be written in a common chart away from the singular locus. 
The lemmas imply that its components converge to the components of $\nabla_u u$ in the limit space.
By assumption $u$ is parallel in the limit space; in particular $\nabla_uu=0$.

This observation will be used to control the term $\<\nabla\phi_n,\nabla_{u_n} {u_n}\>$ at the points far from the singular locus of $\Cone(\theta)\times \RR$.
To do this we only need to assume that $|\nabla\phi_n|$ is bounded.
Next, we describe how to control it near the singularity.

Since $|u_n|=1$, we have $\nabla_{u_n} u_n\perp \nabla f_n$.
Therefore if $\nabla \phi_n$ is proportional to $\nabla f_n$ at some point, then  $\<\nabla\phi_n,\nabla_{u_n} {u_n}\>=0$ at this point.
This observation makes it possible to choose $\phi_n$ so that the term $\<\nabla\phi_n,\nabla_{u_n} {u_n}\>$ vanish around the singular locus of $\Cone(\theta)\times \RR$.
Namely, in addition to the above conditions on $\phi_n$ we have to assume that the identity $\phi_n=\psi\circ f_n$ holds at the points of $M_n$ that are sufficiently close to the singular locus of $\Cone(\theta)\times \RR$.

\begin{wrapfigure}{r}{40 mm}
\vskip-3mm
\centering
\includegraphics{mppics/pic-50}
\vskip0mm
\end{wrapfigure}

Finally, observe that the needed sequence exists.
Indeed, one can take 
\[\phi_n=(\sigma\circ \dist_{p_n})\cdot (\psi\circ f_n)\]
for appropriately chosen fixed mollifiers $\sigma,\psi\:\RR\to \RR$ and $M_n\zz\ni p_n\to p$.

\parit{(\ref{cor:Ricci:vw}).}
Since $G_n\rightharpoonup 0$, we get that the curvature measure of level sets of $f_n$ weakly converges to the curvature of $\Cone(\theta)$.
It follows that 
\[\Qm(u_n,u_n)\rightharpoonup \omega.\]

Suppose $v'_n\perp u_n$ for all $n$.
Part \textit{(\ref{cor:Ricci:Ricci})} implies that  $\Qm(v'_n,v'_n)\zz\rightharpoonup 0$.

Fix $t\in \RR$.
Since the lower bound on sectional curvature of $M_n$ converges to $0$, any partial weak limit of $\Qm(v'_n\zz+t\cdot w_n,v'_n\zz+t\cdot w_n)$ is nonnegative.
It follows that 
\[\Qm(v'_n, w_n)\rightharpoonup 0\] for any sequence of fields $v'_n,w_n$ such that $v'_n\perp u_n$.

Consider the vector fields $v_n',w_n'$ such that 
\begin{align*}
v'_n&\perp u_n,
&
v_n&=a_n\cdot u_n+v_n',
\\
w'_n&\perp u_n,
&
w_n&=b_n\cdot u_n+w_n'.
\end{align*}
Since $\Qm$ is bilinear, we get that
\[\Qm(v_n,w_n)=\Qm(v_n',w_n) + a_n\cdot [\Qm(u_n,w_n')+ b_n\cdot \Qm(u_n,u_n)].\]
By assumption, $a_n=\langle u_n,v_n\rangle $ and $b_n=\langle u_n,w_n\rangle$ uniformly converge to $a$ and $b$ respectively.
Whence the statement follows.
\qeds

\subsection{Three-dimensional case}\label{subsec:3d}

\parit{Proof of the 3-dimensional case in \ref{prop:3parts:codim2+}.}
Suppose $M_n\smooths{}A$ and $\dim A=3$.
Choose a set $Q\subset A$ that admits a bi-Lipschitz embedding into $\RR$. 

Let us split $\mathfrak m$ into negative and positive parts $\mathfrak m=\mathfrak m^+-\mathfrak m^-$; that is,
\[\mathfrak m^\pm(X)\df\sup\set{\pm\mathfrak m(Y)}{Y\subset X}.\]
Since the sectional curvature of $M_n$ is bounded below, we get that $\mathfrak m^-$ has bounded density; in other words, $\mathfrak m^-$ is a regular measure with respect to $\vol^3$ on $A$.
Since $Q$ has zero volume, we get $\mathfrak m^-(Q)=0$.

Set $\mathfrak n=\mathfrak m|_{Q}$; from above we have $\mathfrak n\ge 0$.
By \cite{petrunin-SC}, $\mathfrak n$ is regular with respect to $\vol^1$ on $Q$.
Therefore it is sufficient to show that 
\[(2\cdot\pi-\theta)\cdot (\jac(\bm{h}|_Q))^2\]
is the $\vol^1$-density of $\mathfrak n$
at $\vol^1$-almost all $p\in Q$.

Choose a bi-Lipschitz embedding $s\:Q\to \RR$;
set $K\zz=s(Q)$.
Since $s^{-1}$ and $h\circ s^{-1}$ are Lipschitz, 
by Rademacher's theorem, we can assume that $s^{-1}$ and $h\circ s^{-1}$ are differentiable at almost all $x\in K$.
Moreover, we can assume that $d_xs^{-1}(y)=(\lambda\cdot y,0)\in \RR\times \Cone(\theta)=\T_p$ and the $\vol^1$-density of $\mathfrak n$ at $p=s^{-1}(x)$ is defined.

Shifting and scaling the interval $K$, we may assume that $x=0$ and $\lambda=1$.
In this case, $|\jac_p(h|_Q)|=|d_0(h\circ s^{-1})|$.

Note that we can choose a sequence of points $p_n\in M_n$ and a sequence of factors $c_n\to \infty$ such that $(c_n\cdot M_n,p_n)$ converges to the tangent space~$\T_p=\RR\times \Cone(\theta)$.

Applying Perelman's construction \cite[7.1.1 and 7.2.3]{petrunin-conc} for a horizontal vector in $\T_p=\RR\zz\times \Cone(\theta)$,
we can choose a sequence of functions $f_n\:c_n\cdot M_n\to \RR$ satisfying the assumptions in~\ref{thm:HG-converge}; 
let $u_n=\nabla f_n/|\nabla f_n|$.
Consider the sequence of functions $\hat h_n\:c_n\cdot M_n\to \RR$ defined by 
\[\hat h_n(x)=c_n\cdot(h_n(x)-h_n(p_n)).\]
Since $\lambda=1$, we have that $|\langle\nabla \hat h_n,u_n\rangle|$ uniformly converges to $|d_0(h\zz\circ s^{-1})|$.
By \ref{cor:Ricci}\textit{(\ref{cor:Ricci:vw})}, the sequence of measures $\qm(\hat h_n, \hat h_n)$ on $c_n\cdot M_n\to \RR$ weakly converges to $|d_0(h\zz\circ s^{-1})|^2\cdot\omega_{\Cone(\theta)\times\RR}$.
Recall that $\vol^1$-density of $\omega_{\Cone(\theta)\times\RR}$ on the singular line is $2\cdot \pi-\theta(p)$.

Observe that 
\[\qm(\hat h_n, \hat h_n)[B(p_n,1)_{c_n\cdot M_n}]
=
c_n\cdot \qm(h_n, h_n)[B(p_n,\tfrac1{c_n})_{M_n}]\]
Whence $(2\cdot \pi-\theta(p))\cdot(\jac_p(h|_Q))^2$ is the $\vol^1$-density of $\mathfrak n$ at $p$ as required.
\qeds

\subsection{Higher-dimensional case}\label{subsec:4d}

Suppose $M_n\smooths{} \Cone(\theta)\times \RR^{m-2}$, where $\Cone(\theta)$ denotes a two-dimensional cone with the total angle $\theta<2\cdot\pi$.
First, we will show that the curvatures of $M_n$ in the vertical sectional directions of $\Cone(\theta)\times \RR^{m-2}$ weakly converge to zero;
an exact statement is given in the following proposition.
By combining this result with the 3-dimensional case we get \ref{prop:3parts:codim2+} in all dimensions.

For $X,Y\in \T_p$, denote by $K(X\wedge Y)$ the curvature in the sectional direction of $X\wedge Y$.
A function $f\:\Cone(\theta)\times\RR^{m-2}\to\RR$ will be called a \emph{vertical affine function} if $f$ can be obtained as a composition of the projection to $\RR^{m-2}$ and an affine function on $\RR^{m-2}$.

\begin{thm}{Proposition}\label{prop:vert-vert}
Let $M_n\smooths{} \Cone(\theta)\times\RR^{m-2}$
and $f_n,h_n\: M_n\zz\to \RR$ be sequences of strongly concave smooth Lipschitz functions.
Suppose that $\sec M_n\ge -\tfrac1n$ for each $n$,
and we have pointwise convergences $f_n\to f$ and $h_n\to h$, where $f$ and $h$ are vertical affine functions on $\RR^{m-2}\times \Cone(\theta)$.
Then 
\[K(\nabla f_n\wedge\nabla h_n)\cdot\vol^m \rightharpoonup0.\]

\end{thm}

Let $\Sigma$ be a convex hypersurface in an $m$-dimensional Riemannian manifold $M$.
Suppose $x$ is a smooth point of $\Sigma$; that is, the tangent hyperplane $H_x$ of $\Sigma$ is defined at $x$;
denote by $e_1,\dots,e_{m-1}$ an orthonormal basis of $H_x$.
Set 
\[\Zc_\Sigma(x)= \sum_{i,j} K(e_i\wedge e_j).\]
In other words, 
\[\Zc_\Sigma=\Sc-2\cdot\Ric(n_\Sigma,n_\Sigma),\]
where $n_\Sigma$ is the unit normal vector to $\Sigma$.

Since tangent hyperplanes are defined at almost all points of convex hypersurfaces,
 $\Zc_\Sigma(x)$ is defined almost everywhere on $\Sigma$.

\begin{thm}{Lemma}\label{lem:nonsmooth-convex}
Let $\Sigma$ be a strongly convex hypersurface in an $m$-dimensional Riemannian manifold $M$ with curvature $\ge -1$.
Suppose that for some point $p\in \Sigma$ and $r<1$ the closed ball $\bar B(p,2\cdot r)$ in the intrinsic metric of $\Sigma$ is compact.

Then, 
\[\int_{x\in B(p,r)}\Zc_\Sigma(x)\cdot \vol^{m-1}\le (m-1)\cdot(m-2)\cdot\const(m-1)\cdot r^{m-3},\]
where $\const(m-1)$ is the constant in \ref{cor:Kmax}.
\end{thm}

\parit{Proof.} If $\Sigma$ is smooth, then the inequality follows from \ref{cor:Kmax}, and the fact that curvature cannot decrease when we pass to a convex hypersurface.

In the general case, the surface $\Sigma$ can be approximated by a smooth convex surface;
this can be done by applying the Green--Wu smoothing procedure; compare to \cite{AKP-buyalo}.

It remains to pass to the limit.
More precisely, suppose $\Sigma_n$ is a sequence of smooth  strongly convex hypersurfaces that converges to $\Sigma$.
For each $n$, choose a point $p_n\in\Sigma_n$ such that $p_n\to p$. 
By \cite[Theorem 1.2]{petrunin-QG}, $(\Sigma_n,p_n)$ converges to $(\Sigma,p)$ as in the pointed Gromov--Hausdorff convergence.

Recall that $\vol_{m-1}$ on $\Sigma_n$ weakly converges to $\vol_{m-1}$ on $\Sigma$ (see \cite[10.8]{BGP}).
Further, since $M$ is smooth, $\Zc_{\Sigma_n}$ is bounded in $B(p_n,r)_{\Sigma_n}$.
Therefore, 
\[\int_{x\in B(p_n,r)_{\Sigma_n}}\Zc_{\Sigma_n}(x)\cdot \vol^{m-1}\to \int_{x\in B(p,r)_{\Sigma}}\Zc_{\Sigma}(x)\cdot \vol^{m-1}\qquad\text{as}\qquad n\to\infty\]
follows if for almost all $x\in \Sigma$ we have that
\textit{for any $\eps>0$ there is $\delta>0$ such that
if $x_n\in \Sigma_n$ and $|x_n-x|<\delta$ for large $n$,
then}
\[|\Zc_{\Sigma_n}(x_n)-\Zc_{\Sigma}(x)|<\eps.\]
This condition holds if the tangent plane $H_x$ is defined.
Whence the nonsmooth case follows.
\qeds

\parit{Proof of \ref{prop:vert-vert}.}
Passing to a subsequence, we can assume that
\[K(\nabla f_n\wedge\nabla h_n)\cdot\vol^m \rightharpoonup\mathfrak m\]
for some measure $\mathfrak m$ on $\RR^{m-2}\times \Cone(\theta)$.

First, let us show that $\mathfrak m$ is supported on the singular locus.
If $p$ is not singular, then it has a flat neighborhood.
Therefore by a local version of \ref{A^0} (see also Section~\ref{sec:local})
we get that $\mathfrak m$ vanishes in a neighborhood $U\ni p$.
Indeed, we can include copies of $U$ (which is flat) in the approximating sequence $U_n\subset M_n$ of $U$ and argue as in the introduction.

Let $p$ be a singular point on $\Cone(\theta)\times \RR^{m-2}$;
let us denote its liftings by $p_n\in M_n$.
We can assume that $p$ corresponds to the origin of $\RR^{m-2}$.
Choose points $a_{1,n},\zz\dots,a_{m-2,n}$, $b_{1,n},\zz\dots,b_{m-2,n}$ in $M_n$ such that the functions $f_{i,n}\zz=\dist_{a_{i,n}}\zz-|a_{i,n}-p|$ and $-h_{i,n}\zz=-\dist_{b_{i,n}}+|b_{i,n}-p|$ converge to $i^{\text{th}}$ vertical coordinate function on $\Cone(\theta)\zz\times \RR^{m-2}$.
Further, choose points $c_{1,n},c_{2,n},c_{3,n}$ so that the functions $g_{i,n}\zz=\dist_{c_{i,n}}-|c_{i,n}-p|$ converge to Busemann functions for different horizontal rays in $\Cone(\theta)\times \RR^{m-2}$ emerging from $p$.
Note that the latter implies that the angles $\angk{p_n}{c_{i,n}}{c_{j,n}}$ are bounded away from zero for all large~$n$.

\begin{wrapfigure}{r}{35 mm}
\vskip-0mm
\centering
\includegraphics{mppics/pic-100}
\vskip0mm
\end{wrapfigure}

By \cite[Lemma 7.2.1]{petrunin-conc}, there is an increasing concave function $\phi$ defined in a neighborhood of zero in $\RR$ such that $\phi'$ is close to 1 and for any $\eps>0$ and $i\ne j$ the function 
\[s_{ij,n}
=
\phi\circ g_{i,n}
+
\phi\circ g_{j,n}
+
\sum_i\left[\phi\left(\eps\cdot f_{i,n}\right)+\phi\left(\eps\cdot h_{i,n}\right)\right]\]
is strongly concave in $B(p_n,R)$ for fixed $R>0$ and every large $n$.

Denote by $s_{ij}\:\Cone(\theta)\zz\times \RR^{m-2}\zz\to \RR$ the limits of $s_{ij,n}$.
Note that given $w>0$, we can take small $\eps>0$ so that for all $i\ne j$ the set $s_{ij}^{-1}[-w,w]$ covers the singular locus in $B(p,R)$. 

Note that we can choose $\eps_0>0$ so that for almost all points $x\zz\in B(p_n,R)$ 
the differential $d_xs_{ij,n}$ is linear and $|d_xs_{ij,n}|>\eps_0>0$ for some $i$ and $j$.
Indeed, these differentials are linear outside cutlocuses of $a_{i,n}$, $b_{i,n}$, and $c_{i,n}$;
in particular, they are linear at almost any point $x_n\zz\in M_n$.
Further, if the differential $d_{x_n}s_{12,n}$ is very close to zero,
then the directions of $[x_n,c_{1,n}]$ and $[x_n,c_{2,n}]$ are nearly opposite.
Since $x_n$ is close to $p_n$, we get that for large $n$ the angles $\measuredangle\hinge{x_n}{c_{i,n}}{c_{j,n}}$ is bounded away from zero, we get $|d_{x_n}s_{13,n}|$ is bounded away from zero as well.

Since $f_n$ and $h_n$ are converging to vertical affine functions, we get that for large $n$ their gradients are nearly orthogonal to $\nabla_{x_n} g_{i,n}$ at almost all $x_n\in M_n$.
Suppose $d_xs_{ij,n}$ is linear and $|d_xs_{ij,n}|>\eps_0>0$.
Then gradients $\nabla_{x_n}f_n$ and $\nabla_{x_n}h_n$ are nearly orthogonal to $\nabla_x s_{ij,n}$.

Set $\Sigma_{ij,n}=\Sigma_{ij,n}(c)=\set{x\in M_n}{s_{ij,n}(x)=c}$.
The argument above implies that for almost all points $x\in B(p_n,R)$ one of the sectional directions of the tangent directions $\sigma$ of $\Sigma_{ij,n}$ is close to the sectional direction $\nabla f_n\wedge\nabla h_n$.
In particular, given $\delta>0$, we have 
\[K(\nabla f_n\wedge\nabla h_n)(x)\le K(\sigma)+\delta\cdot|K_{max}(x)|\]
for all large $n$ ($K_{max}$ is defined in \ref{cor:Kmax}).

By \ref{lem:nonsmooth-convex} and the coarea formula, the sum of integral curvatures of $M_n$ in the directions of $\Sigma_{ij,n}$ at $x_n$ at the subsets where $|d_{x_n}s_{ij,n}|>\eps_0$ is bounded by $\const\cdot w$.
By \ref{cor:Kmax}, the same holds for the integral of $K(\nabla f_n\zz\wedge\nabla h_n)$ if $n$ is large.
The proposition follows since $w$ can be taken arbitrarily small.
\qeds

\parit{Proof of the general case of \ref{prop:3parts:codim2+}.}
Choose $m-2$ sequences of strongly concave functions $f_{1,n},\zz\dots, f_{m-2,n}\: M_n\to \RR$ that converge to vertical affine functions $f_1,\zz\dots,f_{m-2}$ on  $\Cone(\theta)\times \RR^{m-2}$ with orthonormal gradients.
It is done using Perelman's construction \cite[7.1.1 and 7.2.3]{petrunin-conc} for the corresponding vertical vectors in $\Cone(\theta)\times \RR^{m-2}$.

Note that the fields $e_{1,n}=\nabla f_{1,n}$, $\dots$, $e_{m-2,n}=\nabla f_{m-2,n}$ are nearly orthonormal;
in particular, they are linearly independent for all large $n$.
Let us add two fields $e_{m-1,n}$ and $e_{m,n}$  so that $e_{1,n},\dots, e_{m,n}$ form a nearly orthonormal frame in~$M_n$;
that is, $\langle e_{i,n}, e_{j,n}\rangle\to 0$ for $i\ne j$ and $\langle e_{i,n}, e_{i,n}\rangle\to 1$ for any $i$ as $n\to\infty$.

Observe that \ref{prop:vert-vert} implies that $K(e_i\wedge e_j)\cdot\vol\rightharpoonup 0$ if $i,j\zz\le m-2$.

Let us show that $K(e_i\wedge e_j)\cdot\vol\rightharpoonup 0$ for $i\le m-2$ and $j\zz\ge m\zz-1$.
The $3$-dimensional case is done already; it is used as a base of induction.
Let us apply the induction hypothesis to the level surfaces of $f_{1,n}$.
(Formally speaking, we apply the local version of the induction hypothesis described in Section~\ref{sec:local}.)
Since the curvature of convex hypersurfaces is larger than the curvature of the ambient manifold in the same direction, we get the statement for $i\ne 1$.
Applying the same argument for the level surfaces of $f_{2,n}$, we get the claim.

Now let us show that $K(e_{m-1}\wedge e_m)\cdot\vol^m\rightharpoonup \omega_{\Cone(\theta)\times \RR^{m-2}}$.
Consider the level sets $L_n$ defined by 
\[f_{1,n}=c_1,\ \dots,\ f_{m-2,n}=c_{m-2}.\eqlbl{eq:fi=ci}\]
Note that $L_n\smooths{} \Cone(\theta)$.
Applying the 2-dimensional case to $L_n$ and the coarea formula, we get that curvatures of $L_n$ weakly converge to $\omega_{\Cone(\theta)\times \RR^{m-2}}$.
It remains to show that the extra term in the Gauss formula for the curvature of $L_n$ weakly converges to zero; 
in other words, the difference between the curvature of $L_n$ and sectional curvature of $M_n$ in the direction tangent to $L_n$ weakly converges to zero.

The 3-dimensional case is done already.
To prove the general case, we apply the 3-dimensional case to the 3-dimensional level sets defined by $m-3$ equations from the $m-2$ equations in \ref{eq:fi=ci}.
(The same argument is used in the proof of \ref{A^0}, and it is written with more details.)

Note that for $\theta=0$, the last argument implies the following:

\begin{thm}{Claim}
Let $M_n\smooths{} A$.
If $A$ has a flat open set $U$,
then $|K_{max}|\zz\cdot\vol_n\zz\rightharpoonup0$ on $U$. 
\end{thm}

In particular, the weak limit of dual curvature tensor has support on the singularity of $\Cone(\theta)\times \RR^{m-2}$.

The same argument as in \ref{cor:Ricci} shows that $\langle\Rm(e_i,e_j)e_q,e_r\rangle\zz\cdot\vol^m\zz\rightharpoonup0$ if at least one of the indices $i,j,q,r$ is at most $m-2$.
The latter statement implies the result.
\qeds

\section{Regular points}\label{sec:ref}

\subsection{Common chart and delta-convergence}\label{subsec:chart+delta}

Choose a smoothing $M_n\smooths{} A$ of an $m$-dimensional Alexandrov space $A$.
Let $p\in A$ be a point of rank $m$; that is, there are $m+1$ points $a_0,\dots, a_m\in A$ such that 
$\angk p{a_i}{a_j}>\tfrac\pi2$ for all $i\ne j$.

Recall \cite[Sec. 7]{petrunin-conc} that we can choose small $r>0$,
finite set of points $\bm{a}_i$ near $a_i$,
and a smooth concave increasing real-to-real function $\phi$ defined on an open interval such that 
\[f_i=\sum_{x\in \bm{a}_i}\phi\circ\widetilde{\dist}_{x,r}\]
is strongly a concave function that is defined in a neighborhood $U\ni p$;
it will be called \emph{smoothed distance chart}.

Since $r$ is small, and all points in $\bm{a}_i$ are near $a_i$ we get that the functions $f_0, \dots, f_m$ are tight in $U$; see the definition in \cite{petrunin-conc}.
In particular, the map $U\to\RR^m$ defined by $x\mapsto (f_1(x),\dots,f_m(x))$ is a coordinate system in $U$.

The presented construction can be lifted to $M_n$.
As a result, we obtain a chart of an open set $U_n\subset M_n$.
Passing to smaller sets we may assume that $U$ and each $U_n$ is mapped to a fixed open set $\Omega\subset\RR^m$ for all large $n$.
Further, we assume that it holds for all $n$; it could be achieved by cutting off the beginning of the sequence $M_n$.

The obtained collection of charts $\bm{x}_n\:U_n\zz\to \Omega$ and $\bm{x}\:U\zz\to \Omega$ will be called a \emph{common chart} at $p$.
It will be used to identify points of $\Omega$, $M_n$, and $A$; 
in addition, we will use it to identify the tangent spaces $\T M_n$ and $\T A$ with $\RR^m$.
For example, we will use the same notation for function $M_n\to \RR$ and its composition $\Omega\to \RR$ with the inverse of the chart $U_n\to \Omega$.
We will use index $n$ or skip it to indicate that the calculations are performed in $M_n$ or $A$ respectively.
For example, given a function $f\:\Omega\to \RR$, we denote by $\nabla_nf$ and $\nabla f$ the gradients of $f\circ \bm{x}_n$ in $M_n$ and $f\circ \bm{x}$ in $A$ respectively.

Recall that $A^\delta$ denotes the set of $\delta$-strained points in $A$.
For a fixed common chart $\bm{x}$ we will use the notation $A^\delta_\Omega$ for the image $\bm{x}(A^\delta)\subset \Omega$.

Part \ref{prop:3parts}\textit{(\ref{prop:3parts:reg})} will follow from certain estimates in one common chart.

\begin{thm}{Definitions}\label{def:delta-converge}
Let $M_n\smooths{} A$, $\dim A=m$;
choose a common chart with range $\Omega\subset \RR^m$.

A sequence of measures $\mathfrak n_n$ defined on $\Omega$ is called \emph{weakly delta-converging}
if the following conditions hold:
\begin{itemize}
 \item Every subsequence of $\mathfrak n_n$ has a weak partial limit.
 \item For any $\eps>0$ there is $\delta>0$ such that for any two weak partial limits $\mathfrak m_1$ and $\mathfrak m_2$ of $(\mathfrak n_n)$ we have  
\[|(\mathfrak m_1-\mathfrak m_2)(S)|<\eps\]
for any Borel set $S\subset A^\delta_\Omega$.
\end{itemize}
A sequence of bounded functions $f_n$ on $\Omega$ is called \emph{weakly delta-converging} if the measures $f_n\cdot \vol_n$ are weakly delta-converging.

A sequence of functions $f_n$ defined on $\Omega$ is called \emph{uniformly delta-converging}
if the following conditions hold:
\begin{itemize}
 \item For any $\eps>0$ there is $\delta>0$ such that such that 
\[\limsup_{n\to\infty} \{f_n(x)\}-\liminf_{n\to\infty}\{f_n(x)\}<  \eps\]
for any $x\in A^\delta_\Omega$.
\end{itemize}

\end{thm}

\begin{thm}{Observation}\label{obs:delta-weak-uniform}
If $f_n$ is uniformly delta-converging and $\mathfrak n_n$ is weakly delta-converging,
then $f_n\cdot \mathfrak n_n$ is weakly delta-converging.
\end{thm}

\subsection{Convergences}\label{sec:test-convergence}

\begin{thm}{Lemma}\label{lem:test-delta}
Let $M_n\smooths{} A$, $\dim A=m$;
choose a common chart with range $\Omega\subset \RR^m$.
Let $f_n\: M_n\to\RR$ be a sequence of $C^1$-functions such that $f_n\zz\smoothto f\:A\to \RR$.
Let us denote by $\partial_1,\dots,\partial_m$ the partial derivatives on $\Omega\subset \RR^m$.
Denote by $g_{ij,n}$ and $g^{ij}_n$ the components of the metric tensors on $M_n$.
Then 
\begin{enumerate}[(i)]
\item\label{lem:test-delta-f} $f_n$ uniformly converges to $f$ on $\Omega$;
\item\label{lem:test-delta-partial} $\partial_if_n$ are uniformly delta-converging;
\item\label{lem:test-delta-g}  $g_{ij,n}$ and $g^{ij}_n$ are uniformly delta-converging for all $i,j$;
moreover, $\det g_{ij,n}$ is bounded away from zero;
\item\label{lem:test-delta|nabla|} $|\nabla_n f_n|$ uniformly delta-converges on $\Omega$; 
\end{enumerate}

\end{thm}

\parit{Proof.}
Part \textit{(\ref{lem:test-delta-f})} is trivial.

\parit{(\ref{lem:test-delta-partial}).}
Suppose $a_0, a_1,\dots,a_m$ struts $p$ (see the definition in \cite{AKP}), and the geodesics $[pa_i]$ are uniquely defined.
In this case, for any sequence of points $a_{i,n}, p_n\in M_n$ such that $a_{i,n}\to a_i$, and $p_n\to p$ as $n\to\infty$, we have
\[\lim_{n\to\infty}\measuredangle \hinge{p_n}{a_{i,n}}{a_{j,n}}\ge \measuredangle \hinge{p}{a_i}{a_j}.\]

If $\T_p$ is Euclidean, then ($n$+1)-point comparison implies that equality holds in the last inequality.

Note that the angles $\measuredangle \hinge{p_n}{a_{i,n}}{a_{j,n}}$ for all $i,j>0$ completely describe the metric tensor at $p_n$ in the basis $V_{1,n},\dots,V_{m,n}$, where $V_{i,n}$ is the unit vector in the direction of $[p_n,a_{i,n}]$.

If $f_n\zz\smoothto f$, then $V_{i,n}f_n$ completely describes $\nabla_{p_n}f_n$ in the basis $V_{1,n},\zz\dots,V_{m,n}$.
From above, we can express $|\nabla_{p_n}f_n|$ in terms of $V_{i,n}f_n$ and the angles $\measuredangle \hinge{p_n}{a_{i,n}}{a_{j,n}}$.
Whence we get convergence $|\nabla_{p_n} f_n|\zz\to |\nabla_{p} f|$ and therefore 
\[\langle \nabla_{p_n} f_n,\nabla_{p_n} h_n\rangle\zz\to \langle \nabla_p f,\nabla_p h\rangle\] if $h_n\smoothto h$; the latter follows by identity $4\cdot B(x,y)\zz=B(x\zz+y,x\zz+y)\zz-B(x\zz-y,x\zz-y)$ for any bilinear form $B$.

Note that the partial derivatives $\partial_if_n$ at a regular point $p$ can be expressed in terms of $\langle d_pf_n,d_px_j\rangle_n$ and $\langle d_px_j,d_px_k\rangle_n$, where $x_1,\dots,x_m$ are the coordinate functions of the chart.
Therefore, we get that all $\partial_if_n$ converge at any regular point.

Finally, observe that if $p$ is a $\delta$-strained point for sufficiently small $\delta>0$,
then the calculations above go thru with a small error.
Whence the statement follows.

\parit{(\ref{lem:test-delta-g}).}
This part follows from the proof of \textit{(\ref{lem:test-delta-partial})} since $g^{ij}_n\zz=\langle d_px_i,d_px_j\rangle_n$ and $g_{ij,n}$ can be expressed thru~$g^{ij}_n$.

\parit{(\ref{lem:test-delta|nabla|}).}
Note that $|\nabla_n f_n|$ can be expressed from $g^{ij}_n$ and $\partial_if_n$.
Since these quantities are delta-converging, so is $|\nabla_n f_n|$.
\qeds

The following lemma relies on the DC-calculus which is discussed in Section~\ref{sec:DC};
this section includes the definition of DC and DC$_0$ functions, as well as DC convergence.
Since \textit{test convergence implies DC convergence} (see \ref{obs:test-DC}),
the lemma also holds for test-converging sequences of functions.

\begin{thm}{Lemma}\label{lem:test-delta-partial-g}
Let $M_n\smooths{} A$.
Choose a common chart $\bm{x}_n\:U_n\subset M_n\zz\to \Omega$ and $\bm{x}\:U\subset A\zz\to \Omega$ with range $\Omega\subset \RR^m$.
Let $f_n\: M_n\to\RR$ be a sequence of smooth functions that DC converges to a DC$_0$ function $f\:A\to \RR$.
Let us denote by $\partial_1,\dots,\partial_m$ the partial derivatives on $\Omega\subset \RR^m$.
Denote by $g_{ij,n}$ and $g^{ij}_n$ the components of the metric tensors on $M_n$.
Then  the partial derivatives $\partial_kg_{ij,n}$, $\partial_k g^{ij}_n$,  $\partial_j\partial_if_n$, as well as their products to uniformly delta-converging functions,  are weakly converging.

\end{thm}

\parit{Proof.} 
The weak convergence of $\partial_kg_{ij,n}$, $\partial_k g^{ij}_n$, and $\partial_j\partial_if_n$
follows from \ref{metricBV}.
By Observation~\ref{obs:delta-weak-uniform}, products of these partial derivatives to uniformly delta-converging sequences of functions are weakly delta-converging.

Let $h_n\:M_n\to\RR$ be a uniformly delta-converging sequence.
Note that its limit is well defined in $A^\circ$; denote it by $h$;
let us extend it by $0$ to the whole $A$.

Denote by $\mathfrak m_n$ one of the measures on $\Omega$ with the density $\partial_kg_{ij,n}$, $\partial_k g^{ij}_n$, or $\partial_j\partial_if_n$.
Let $\mathfrak m$ be the corresponding limit measure $\partial_kg_{ij}$, $\partial_k g^{ij}$, or $\partial_j\partial_if$.
We need to show that 
\[\int_{\Omega} (h_n\circ\bm{x}_n^{-1}) \cdot \phi\cdot \mathfrak m_n
\to
\int_{\Omega} (h\circ\bm{x}^{-1}) \cdot \phi\cdot \mathfrak m\quad\text{as}\quad n\to\infty
\eqlbl{eq:convegence}\]
for any continuous function $\phi\:\Omega\to\RR$ with compact support.

Choose $\eps>0$; let $\delta>0$ be as in \ref{def:delta-converge} (for $h_n$).
The set $S^\delta_\Omega=\Omega\setminus A^\delta_\Omega$ is a closed subset of $\Omega$.
By~\ref{metricBV+}, $|\mathfrak m|(S^\delta_\Omega)=0$.
Therefore we can choose an open neighborhood $N\subset \Omega$ of $S^\delta_\Omega$ such that $|\mathfrak m|(N)<\eps$.
Choose two nonnegative continuous functions $\phi_0$ and $\phi_1$ such that 
\[\phi=\phi_0+\phi_1,
\qquad
\supp \phi_0\subset N,
\qquad
\supp \phi_1\subset A^\delta_\Omega=\Omega\setminus S^\delta_\Omega.\]

Note that the sequence $a_n=\int_{\Omega} (h_n\circ\bm{x}_n^{-1}) \cdot \phi_0\cdot \mathfrak m_n$ \textit{converges with error} 
$\eps_0=\eps\cdot c \cdot\max\{\,|\phi|\,\}$, where $c$ is a bound on $|h_n|$. 
In other words, the upper and lower limits of $a_n$ differ by at most~$\eps_0$.
Similarly, $b_n=\int_{\Omega} (h_n\circ\bm{x}_n^{-1}) \cdot \phi_1\cdot \mathfrak m_n$ converges with error $\eps_1=\eps\cdot |\mathfrak m|\cdot c \cdot\max\{\,|\phi|\,\}$.
Since $\eps>0$ is arbitrary, we get \ref{eq:convegence}.
\qeds

\subsection{Proof modulo key lemma}\label{subsec:modulo-key}

\begin{thm}{Key lemma}\label{A^0}
Choose a common chart with range $\Omega\subset \RR^m$ for a smoothing $M_n\smooths{} A$.
Choose a component $\Rm_{ijsr,n}$ of the curvature tensor of $M_n$ in $\Omega$.
Then $\Rm_{ijsr,n}\cdot \vol^m_n$ is a weakly delta-converging sequence of measures.
\end{thm}

The proof of the key lemma will take the remaining part of this section;
in the current subsection, we show that it implies \ref{prop:3parts}\textit{(\ref{prop:3parts:reg})}.

\parit{Proof of \ref{prop:3parts}(\ref{prop:3parts:reg}) modulo \ref{A^0}.}
Recall that components of $\qm_n$ can be expressed from the components of $\Rm_n$.
Therefore, the key lemma implies delta-convergence of components of $\qm_n$.

Choose sequences of test functions $f_{1,n},\zz\dots,f_{m-2,n},h_{1,n},\zz\dots,h_{m-2,n}$ on $M_n$ that test-converge to $f_{1},\zz\dots,f_{m-2},h_{1},\zz\dots,h_{m-2}\:A\to \RR$.
By \ref{lem:test-delta}, we have delta-convergence of the partial derivatives $\partial_if_{j,n}$ and 
$\partial_ih_{j,n}$ to $\partial_if_{j}$ and 
$\partial_ih_{j}$ respectively.
The measures $\qm_n(f_{1,n},\zz\dots,f_{m-2,n},\zz{}h_{1,n},\zz\dots,h_{m-2,n})$ 
can be expressed as a linear combination of the components of $\qm_n$ with coefficients expressed in terms of $\partial_if_{j,n}$.
By \ref{obs:delta-weak-uniform}, it follows that the sequence of measures 
\[\mathfrak m_n=\qm_n(f_{1,n},\dots,f_{m-2,n},h_{1,n},\zz\dots,h_{m-2,n})\]
is delta-converging.

Finally, recall that 
\[A^\circ=\bigcap_{\delta>0}A^\delta.\]
Therefore delta-convergence of $\qm_n(f_{1,n},\dots,f_{m-2,n},h_{1,n},\zz\dots,h_{m-2,n})$ implies \ref{prop:3parts}\textit{(\ref{prop:3parts:reg})}.
\qeds

\subsection{Strange curvature}

Suppose $M$ is a 3-dimensional Riemannian manifold.
\emph{Strange curvature tensor} $\Str$ on $M$ is a bilinear form that is uniquely defined by
$$\Str(w,w)=\Sc\cdot |w|^2-\Ric(w,w)$$
for $w\in \T M$.
Note that $\Str$ completely describes the Ricci curvature tensor $\Ric$.
Further, since $M$ is 3-dimensional, $\Str$ completely describes the curvature tensor $\Rm$ of $M$.

In Riemannian manifolds, we can (and will) use the metric tensor to identify tangent and cotangent bundles.
Therefore the tensor $\Str$ can be applied to vector fields and forms;
in particular, for any smooth function $f$ we have
\[\Str(df,df)
=
\Str(\nabla f,\nabla f).\]

\begin{thm} {Proposition}\label{strconvergence}
Let $M_n\smooths{} A$ and $\dim A=3$; choose a common chart with range $\Omega\subset \RR^3$.
Suppose that $f$ is a convex combination of coordinate functions of the chart.
Then the measures 
\[\mathfrak m_n=\Str_n(df,df) \cdot \vol^3_n\] are weakly delta-converging in $\Omega$.

\end{thm}

The definition of strange curvature tensor is motivated by the following integral expression from \ref{thm:bochner-formula}:
$$\int\limits_\Omega \phi\cdot \Str(u, u)
=\int\limits_\Omega \phi\cdot \Int+
\int\limits_\Omega \l[H\cdot\<u,\nabla\phi\>- \<\nabla\phi,\nabla_{u} u\> \r],\eqlbl{eq:Bochner2}$$
where 
\begin{itemize}
\item $u=\nabla f/|\nabla f|$,
\item $H(x)$ --- the mean curvature of the level set $f^{-1}(f(x))$,
\item $\Int(x)$ --- the scalar curvature of  $f^{-1}(f(x))$.
\end{itemize}
This formula is the main tool in the proof of the proposition.
It reduces the proposition to the following two lemmas;
each lemma provides the convergence of an integral term in the right-hand side of \ref{eq:Bochner2}.

\begin{thm} {Lemma}\label{Int}
In the assumptions of Proposition~\ref{strconvergence}, $\Int_n\zz\cdot \vol^3_n$ is a delta-converging sequence of measures on $\Omega$.
\end{thm}

The proof of this lemma uses the convergence of curvature measures
$\Int_n\cdot \vol^2$
on  the 2-dimensional level sets of concave functions $f$ and the coarea formula.
Recall that the sequence $|\nabla_n f|$ is only weakly delta-converge (see \ref{lem:test-delta}\textit{(\ref{lem:test-delta|nabla|})}).
Since the factor $|\nabla_n f|$ appears in the coarea formula,
we get that only weak delta-convergence of $\Int_n\cdot \vol^3_n$.

\parit{Proof.}
Recall that any point  in an Alexandrov space $A$ has a convex neighborhood \cite{petrunin-conc}.
This construction can be lifted  to the smoothing sequence~$(M_n)$.
Let 
$V\subset A$ be an open  convex
neighborhood of $x$ and
$ V_n\subset M_n$
be open convex sets such that
$\bar V_n  \GHto   \bar V$.

Set
\begin{align*}
L_{t,n}&=f^{-1}(t)\cap V_n,&
C_{t,n}&=f^{-1}[t,\infty)\cap \bar V_n,
\\
L_{t}&=f^{-1}(t)\cap V,&
C_{t}&=f^{-1}[t,\infty)\cap \bar V.
\end{align*}

For every $t$ and $n$, the set $C_{t,n}$ is a convex subset in Alexandrov space 
 and hence is an Alexandrov space 
 with curvature $\ge -1$.
Note that
$C_{t,n} \GHto C_{t}$.
Let us equip the boundaries  $\partial C_{t,n}$ and
 $\partial C_{t}$ with the induced inner metrics.
By \cite[Theorem 1.2]{petrunin-QG}, $\partial C_{t, n}$ converges to $\partial C_{t}$ as $n\to\infty$.

By \cite{AKP-buyalo},
$\partial C_{t,n}$ is
an Alexandrov space 
with curvature $\zz\ge -1$;
hence, so is the limit
$\partial C_{t}$.

Note that $L_{t,n}$ with induced inner metric is isometric to its image in $\partial C_{t,n}$.
Since 
$\partial C_{t}$
is an extremal subset of
$C_{t}$, the inner metric of
$\partial C_{t} $ is bi-Lipschitz to 
the metric restricted from $A$.
It follows that
we can take $r$ sufficiently small
such that for all $t$ and
$U_{t,n}=L_{t,n}\cap B(x_n,r)$
we will have
\[h_n\le \tfrac 1{10}\cdot \dist   (U_{t,n},\partial C_{t,n}\setminus L_{t,n}),\]
where $h_n$ denotes the intrinsic diameter of $U_{t,n}$.
Then the local version of the 2-dimensional case of the main theorem can be applied to $U_{t,n}$; it implies weak convergence of measures  $\Int_n\cdot \vol^2_n$ on $L_{t,n}$.

Choose a smooth function
$\phi\:B(x_n,r)\to\RR$ with a compact support in $A^\delta_\Omega$.
Applying the coarea formula, we get
$$\int\limits_{s\in \Omega}\Int_n(s)\cdot\phi(s)\cdot\vol^3_n
=\int\limits_{-h}^{h} d t\cdot 
\int\limits_{s\in U_{t,n}}
 \frac{\phi_n(s)}{|\nabla_n f(s)|}\cdot\Int_n(s)\cdot \vol^2.\eqlbl{eq:Int-coarea}$$
 
Note that $\nabla_n f$ is bounded away from zero.
By \ref{lem:test-delta}\textit{(\ref{lem:test-delta|nabla|})}, $\frac{1}{|\nabla_n f(s)|}$ is uniformly delta-converging.
Recall that $\Int_n\cdot \vol^2$ are weakly converging measures on $L_n$ \cite[VII \S13]{AZ}.
Therefore \ref{obs:delta-weak-uniform} implies
that $\Int_n\cdot\vol^3_n$ is weakly delta-converging measures on $M_n$.
\qeds

The following lemma is related to the convergence of the second integral in \ref{eq:Bochner2}, the proof uses
the DC calculus in a common chart; see Section~\ref{sec:DC}.

\begin{thm}{Lemma}\label{HnablaU}
In the assumptions of Proposition~\ref{strconvergence}, 
suppose $\phi\:\Omega\zz\to\RR$ is a smooth function with compact support.
Then
\[\int\limits_\Omega \l[H_n\cdot\<u_n,\nabla_n\phi\>_n- \<\nabla_{u_n} u_n,\nabla_n\phi\>_n \r]\cdot\vol^3_n
\eqlbl{eq:HnablaU}\]
converges, where $H_n$ and $u_n$ as in \ref{eq:Bochner2}.
\end{thm}

\parit{Proof.}
Note that $H_n=\div u_n$.
Let us rewrite the first term of \ref{eq:HnablaU} in coordinates:
\begin{align*}
\int\limits_{\Omega}
&\biggl[
\sum_{i}
\Bigl(\partial_i u_n^i
+
\tfrac{u^i_n}{2}\cdot\partial_i \log \det g_{n}\Bigr)\biggr]
\cdot
\biggl[\sum_{i,j}  u_n^i\cdot \partial_i \phi\biggr]
\cdot\sqrt {\det g_n}
\cdot dx^1 dx^2 dx^3.
\end{align*}
We also have
$$u_n^i
=
\frac{\sum_j g^{ij}\cdot\partial_j f}
{\sqrt{\sum_{j,k}g^{jk}\cdot\partial_j f\cdot\partial_k f}}.$$

Taking the derivatives, we see under the integral a sum of products the following two types of expressions:
the first a partial derivative
$\partial_kg_{ij,n}$,
$\partial_kg^{ij}_n$,
or $\partial_i\partial_j f$,
and the second is an expression made from
$g_{ij,n}$,
$g^{ij}_n$,
$\partial_i f$,
$\partial_i \phi$.
Applying \ref{lem:test-delta} and \ref{lem:test-delta-partial-g}, we get that 
the integral converges.

Further, for the second term in \ref{eq:HnablaU}
we have
\begin{align*}
\int\limits_{M_n}&\<\nabla\phi_n,\nabla_{u_n} u_n\>\cdot\vol^3_n=\\
&=
\int\limits_{\Omega}
dx^1 dx^2 dx^3\cdot
\sum_{i,j,k}
u^i_n\cdot \partial_k \phi\cdot\sqrt{\op{det}g_{n}}\cdot
\biggl(
\partial_i u^k_n+
\tfrac{u^j_n}2
\cdot
\sum_s
\left(
\partial_i g_{ js,n}
+
\partial_j g_{si,n}
-
\partial_s g_{i j,n}\right)
\cdot
g^{ks}_n
\biggr).
\end{align*}
The convergence follows by the same argument.
\qeds

\parit{Proof of Proposition~\ref{strconvergence}.}
By \ref{eq:Bochner2}, \ref{Int}, and \ref{HnablaU} we get that $\Str_n(u_n,u_n)\zz\cdot \vol^3_n$ is a weakly delta-converging sequence of measures.
It remains to apply \ref{lem:test-delta}\textit{(\ref{lem:test-delta|nabla|})} and
\ref{obs:delta-weak-uniform}.
\qeds

\subsection{Three-dimensional case}\label{sec:3D-smooth}

In this section, we prove Lemma~\ref{A^0} in the 3-dimensional case.

Vectors $w_1,\dots,w_{m(m+1)/2}\in\R^m$ are said to be \emph{in general position}
if the vectors $w_i\otimes w_i$ form a basis in $\RR^{m\cdot(m+1)/2}$ --- the symmetric square of $\RR^m$.
In this case, any quadratic form $Q$ on $\R^m$
can be computed from the $\tfrac{m(m+1)}2$ values 
\[Q(w_1, w_1),\ \dots,\ Q(w_{m(m+1)/2}, w_{m(m+1)/2}).\] 
More precisely, there are rational functions 
$s_1,\dots,s_{m(m+1)/2}$ that take $\tfrac{m(m+1)}2$ vectors in $\R^m$ and return a quadratic form on $\R^m$
such that
$$Q=\sum_{k=1}^{m(m+1)/2}s_k(w_1,\dots,w_{m(m+1)/2})\cdot Q(w_k,w_k).
\eqlbl{Qij}$$

Note that the vectors $w_1,\zz\dots,w_{m(m+1)/2}\in\R^m$  are  in general position if and only if 
$s_k(w_1,\zz\dots,w_{m(m+1)/2})$ are finite for all $k$.
Since $s_k$ are rational functions, we get the following:

\begin{thm}{Observation}\label{obs:genpos}
Suppose that vectors $w_1,\zz\dots,w_{m(m+1)/2}\zz\in\R^m$ are in general position.
Then the functions $s_1,\zz\dots,s_{m(m+1)/2}$ are Lipschitz in a neighborhood of $(w_1,\zz\dots,w_{m(m+1)/2})\zz\in (\R^m)^{m(m+1)/2}$.
\end{thm}

\parit{Proof of the 3-dimensional case of \ref{A^0}.} 
Choose a common chart
\[M_n\supset U_n\to\Omega,
\quad\text{and}\quad
A\supset U\to\Omega.\]
Let us use it to identify tangent spaces of $M_n$ and $A$ with $\RR^3$.
 
Choose 6 sequences of convex combinations of coordinate functions $f_1,\dots,  f_6$, such that
$\nabla f_1,\zz\dots,\nabla f_6$ are in general position at $p\in \Omega$.
We can assume that $\Omega$ is a small neighborhood of $p$, so by Proposition~\ref{strconvergence}
the measures $\Str_n(\nabla_n f_k, \nabla_n f_k)\zz\cdot\vol^3_n$ weakly delta-converges
on $A^\delta_\Omega$ for $k=1,\dots,6$.

By \ref{obs:genpos}, the functions $s_i$ are Lipschitz in a neighborhood of $(\nabla f_1,\zz\dots,\nabla f_6)\zz\in (\R^3)^6$.
Applying \ref{Qij}, we get that 
\begin{align*}
\Str
&=
\sum_{k=1}^6 s_k(\nabla_n f_1,\dots, \nabla_n f_6)\cdot\Str_n(\nabla_n f_k,\nabla_n f_k).
\end{align*}
Hence the measure $\Str_n(dx_i,dx_j)\cdot\vol^3_n$ are weakly delta-converging for all $i$ and $j$,
where $x_1,x_2,x_3$ is the standard coordinates in $\RR^3$.

By Lemma~\ref{lem:test-delta}, the sequence of metric tensors $g_n$ of $M_n$ on $\Omega$ is uniformly delta-converging.
Since the following equality
\[\operatorname{Tr}\Str_n=\sum_{i,j}g_{ij,n}\cdot\Str_n(dx_{i},dx_{j})\]
holds almost everywhere, we get that the sequence of measures $\operatorname{Tr}\Str_n\zz\cdot\vol^3_n$ is weakly delta-converging.

Note that for $3$-dimensional manifolds we have
$$\Qm_n(V,V)=\Str_n(V,V)-\tfrac14\cdot |V|^2\cdot \Tr \Str_n.\eqlbl{Q}$$
Hence the measures $\Qm_n(d x_i, d x_j)\cdot\vol^3_n$ are weakly delta-converging for all $i$ and~$j$.

Finally, according to \ref{lem:test-delta}\textit{(\ref{lem:test-delta-partial})}, the components $\alpha_{ik,n}$ of $\nabla_n f_k$ are uniformly delta-converging.
The result follows since
\[\qm(f_k, f_k)=
\sum_{i,j} \alpha_{ik,n}\cdot \alpha_{jk,n}\cdot \qm(x_{i}, x_{j}).\]
\qeds

\subsection{Higher-dimensional case}

\begin{thm}{Observation}\label{obs:nested-convex}
Choose a common chart with the range $\Omega\subset\RR^m$ for a smoothing $M_n\smooths{} A$.
Consider the sequence of coordinate level sets $\Omega\zz=L_m\zz\supset L_{m-1}\zz\supset\zz\dots\zz\supset L_0$, 
where $L_i\zz=L_i(c_{i+1},\dots,c_m)$ is defined by setting the last $m-i$ coordinates to be $c_{i+1},\dots,c_m$ respectively.
Then the level sets $L_i$ is a smooth convex hypersurface in $L_{i+1}$ in each $M_n$;
in particular, each $L_i$ has sectional curvature bounded below by $-1$.

Moreover, there is an open set $O$ in the space of linear transformations of $\RR^n$
such that 
the same holds after applying any linear transformation $T\in O$ to $\Omega$.  
\end{thm}

\parit{Proof of the key lemma (\ref{A^0}).}
Let us use notations as in the observation.
By the key lemma (\ref{A^0}) in dimensions $2$ and $3$,  we get weak delta-convergence of curvature tensors on $L_2$ and $L_3$.
(Again, we apply the local version of these statements as described in Section~\ref{sec:local}.)
In particular, applying the coarea formula, we get convergence of sectional curvatures of $L_3$ in the directions of $L_2$ as well as 
the sectional curvature of $L_2$ 
for all values $c_3,\dots,c_m$.
The difference between these curvatures is the Gauss curvature $G_n$ of $L_2$ as a submanifold in $L_3$.
Therefore, $G_n$ is weakly delta-converging as well.

Consider a linear transformation of $\Omega$ that preserves the direction of $L_2$.
By the last statement in \ref{obs:nested-convex},
the above argument shows weak delta-convergence of $G_n(w)$, where the direction $w$ of $L_3$ on $L_2$ can be chosen in an open set of $\RR^{m-2}$ --- the space transversal to $L_2$.
In particular, we may choose directions $w_1,\dots, w_{(m-2)\cdot(m-1)/2}$ in $\RR^{m-2}$ that form a generic set (see the  definition in Subsection~\ref{sec:3D-smooth}).

Denote by $G_n^+$ the term in the Gauss formula for $L_2$ in $M_n$;
that is, $G_n^+$ is the difference between the curvature of $L_2$ and the sectional curvature of $M_n$ in the same direction.
Denote by $g_n$ the Riemannian metric of $M_n$ in $\Omega$.
Note that 
\[G_n^+=\sum\alpha_{k,n}\cdot G_n(w_k),\]
where the coefficients $\alpha_{k,n}$ depend continuously on $w_1,\zz\dots, w_{(m-2)\cdot(m-1)/2}$, and the components of $g_n$. 
It follows that weak delta-convergence of $G_n(w_k)$ implies weak delta-convergence of $G_n^+$ as $n\to\infty$.
Since the curvature of $L_2$ is weakly delta-converging, it implies weak delta-convergence of sectional curvature in the direction of $L_2$.

By the second statement in the observation,
the above argument can be repeated after applying a linear transformation of $\Omega$ that changes the direction of $L_2$ slightly.
It follows that sectional curvatures converge for a generic array of simple bivectors in $\RR^m$.
Note that the curvature tensor can be expressed from these sectional curvatures and the metric tensor.
Hence, the weak delta-convergence of components of curvature tensor and therefore dual curvature tensor follows.
\qeds

\section*{Details}
\addcontentsline{toc}{part}{Details}

\section{Bochner formula}\label{sec:bochner}

Let $M$ be a Riemannian $m$-manifold and $f\:M\to\R$ be a smooth function without critical points on an open domain $\Omega\i M$.
Set $u=\nabla f/|\nabla f|$.
Let us define $\Int_f(x)$ (or just $\Int$) to be scalar curvature of the level set $L_x=f^{-1}(f(x))$ at $x\in L_x\i M$.
Set
\begin{enumerate}
 \item $\kappa_1(x)\le\kappa_2(x)\le\dots\le\kappa_{m-1}(x)$ the principal curvatures of $L_x$ at $x$;
 \item $H=H_f(x)=\kappa_1+\kappa_2+\dots+\kappa_{m-1}$ is mean curvature of $L_x$ at $x$
\item $G=G_f(x)=2\cdot\sum_{i<j}\kappa_i\cdot\kappa_j$ is the extrinsic term
 in the Gauss formula for $\Int_f(x)$. 
\end{enumerate}

Recall that the strange curvature $\Str$ is defined as
\[\Str(u)=\Sc-\<\Ric(u),u\>,\]
where $\Sc$ and $\Ric$ denote scalar and Ricci curvature respectively.

\begin{thm}{Bochner's formula}\label{thm:bochner-formula}
Let $M$ be an $m$-dimensional Riemannian manifold,
 $f\:M\to\R$ be a smooth function without critical points on an open domain $\Omega\subset M$, and $u=\nabla f/|\nabla f|$.
Assume $\phi\:\Omega\to\R$ is a smooth function with compact support.
Then 
\[\int\limits_\Omega \phi\cdot \langle \Ric u,u\rangle =
\int\limits_\Omega [\phi \cdot G+H\cdot\<u,\nabla\phi\>-\<\nabla\phi,\nabla_{u} {u}\>]
\eqlbl{eq:Bochner(-1)}
\]
and
\[\int\limits_\Omega \phi\cdot \Str(u)
=
\int\limits_\Omega \l[H\cdot\<u,\nabla\phi\>- \<\nabla\phi,\nabla_u u\> \r]+
\int\limits_\Omega \phi\cdot \Int_f.
\eqlbl{eq:Bochner0}\]
\end{thm}

The following calculations are based on \cite[Chapter II]{lawson-michelsohn}.
The Dirac operator will be denoted by $D$.
We use the Riemannian metric to identify differential forms and multivector fields on $M$.
Therefore the statement about differential forms can be also formulated in terms of multivector fields and the other way around.

\parit{Proof.}
Assume $b_1,\dots, b_m$ is an orthonormal frame such that $b_m=u$, 
then 
\[\Sc-2\cdot \<\Ric(u),u\>=2\cdot \sum_{i<j<m} \sec(b_i\wedge b_j).\] 
Therefore the Gauss formula can be written as
\[
\begin{aligned}
\Int&=G+\Sc-2\cdot \<\Ric(u),u\>=
\\
&=G+\Str(u)- \<\Ric(u),u\>.
\end{aligned}
\eqlbl{eq:gauss}
\]

We can assume that $b_i(x)$ points in the principal directions of $L_x$ for $i<m$;
so we have $\nabla_{b_i}u=\kappa_i\cdot b_i$ at $x$.
We will denote by ``$\,\bullet \,$'' the Clifford multiplication;
recall that $b_i\bullet  b_i=-1$.
Note that 
\begin{align*}
Du&=\sum_{i} b_i\bullet  \nabla_{b_i}u=
\\
&=\sum_{i<m}\kappa_i\cdot  b_i\bullet  b_i+u\bullet  \nabla_{u}u=
\\
&=-H+u\bullet  \nabla_{u}u.
\end{align*}
Since $\<\nabla_u u,u\>=0$, we get $H\perp (u\bullet  \nabla_{u}u)$.
Therefore
\begin{align*}
\langle Du,Du \rangle&=
\l(\sum_{i<m}\kappa_i\r)^2+|u\bullet  \nabla_{u}u|^2=
\\
&=H^2+|\nabla_{u}u|^2.
\end{align*}
On the other hand
$$\nabla u=\sum_{i<m}\kappa_i\cdot b_i\otimes b_i+\nabla_u u\otimes u,$$
hence
$$\langle\nabla u,\nabla u\rangle =
\sum_{i<m}\kappa_i^2+|\nabla_{u}u|^2.$$
Therefore
$$\langle D u,D u\rangle-\langle \nabla u,\nabla u \rangle =2\cdot\sum_{i<j}\kappa_i\cdot\kappa_j=G.$$

Following the calculations in \cite[II.5.3]{lawson-michelsohn}, we get
\begin{align*}
\int\limits_\Omega\phi\cdot\l[\<D u,D u\>-\<D^2 u, u\>\r]
&=
-\int\limits_\Omega\<\nabla\phi\bullet u,D u\>
=
\\
&=
-\int\limits_\Omega\l[H\cdot\<\nabla\phi,u\>- \<\nabla\phi,\nabla_u u\> \r].
\end{align*}

Since $|u|\equiv 1$, we have $\<\nabla_{\nabla\phi}  u, u\>=0$.
Therefore
$$\int\limits_\Omega\phi\cdot\l[\<\nabla u,\nabla u\>-\<\nabla^*\nabla u, u\>\r]
=
\int\limits_\Omega\<\nabla_{\nabla\phi}  u, u\>=0.$$

By Bochner formula \cite[II.8.3]{lawson-michelsohn},
$$D^2u-\nabla^*\nabla u=\Ric(u);$$
in particular, 
$$\phi\cdot \<D^2u,u\>-\phi\cdot \<\nabla^*\nabla u,u\>=\phi\cdot \<\Ric(u),u\>.\eqlbl{eq:prebochner}$$
Integrating \ref{eq:prebochner} and applying the derived formulas, we get
\begin{align*}
\int\limits_\Omega \phi \cdot G
&=\int\limits_\Omega\phi\cdot\left[\<D u,D u\>-\<\nabla u,\nabla u\>\right]
=
\\
&=\int\limits_\Omega \phi\cdot \Ric(u,u) 
-
\int\limits_\Omega\left[ H\cdot\<u,\nabla\phi\>- \<\nabla\phi,\nabla_u u\>\right].
\end{align*}

It remains to apply the Gauss formula \ref{eq:gauss}.
\qeds

\section{DC-calculus}\label{sec:DC}

Let $f$ be a continuous function defined on an open domain of an $m$-dimensional Alexandrov space~$A$.
Recall that $f$ is \emph{DC} if it can be presented locally as a difference between two concave functions.
Recall that for any point $p\zz\in A$ there is a $(-1)$-concave function defined in a  neighborhood  of~$p$ \cite[3.6]{PerMorse}.
Therefore we can say that \textit{$f$ is DC if and only if it can be presented locally as a difference between two semiconcave functions.} 

Suppose that a sequence of Alexandrov spaces $A_n$ converges to Alexandrov space $A$ without collapse.
Let $f_n$ and $f$ be DC functions defined on open domains $\Dom f_n\zz\subset A_n$ and $\Dom f \zz\subset A$.
Suppose that for any $p\in \Dom f$ there is a sequence $p_n\in \Dom f_n$ and $R>0$ such that $p_n\to p$ and $B(p_n,R)_{A_n}\subset \Dom f_n$, $B(p,R)_{A}\subset\Dom f$
and for some fixed $\lambda\in\RR$, and each large $n$ we have $\lambda$-concave functions $a_n$ and $b_n$ defined in $B(p_n,R)_{A_n}$ and $\lambda$-concave functions $a$ and $b$ defined in $B(p,R)_{A}$
such that $f_n=a_n-b_n$ and $f=a-b$ and the sequences
$a_n$ and $b_n$ converge to functions $a$ and $b$ respectively.
In this case, we say that $f_n$ is \emph{DC-converging} to $f=a-b\:A\to \RR$ as $n\to\infty$; briefly $f_n\DCto f$.

A DC function $f\:A\to \RR$ is called \emph{DC$_0$} if it is continuously differentiable in $A^\circ$.
More preciously, for any smoothed distance chart $\bm{x}\:U\subset A \to \RR^m$ (see Section~\ref{subsec:chart+delta}) the restriction $f\circ\bm{x}^{-1}|_{\bm{x}(A^\circ)}$ is continuously differentiable.

{\sloppy

\begin{thm}{Observation}\label{obs:test-DC}
Any test function is DC$_0$.
Moreover, test convergence implies DC-convergence. 
\end{thm}

}

\parit{Proof.}
Choose a test function $f=\phi(\widetilde\dist_{p_1,r},\dots,\widetilde\dist_{p_n,r})$.
Note that the function $\phi$ can be presented locally as a difference between $C^2$-smooth concave functions increasing in each argument; say $\phi=\psi-\chi$.

For the first part of the observation, it remains to observe that the functions \[a=\psi(\widetilde\dist_{p_1,r},\zz\dots,\widetilde\dist_{p_n,r}),\qquad b=\chi(\widetilde\dist_{p_1,r},\zz\dots,\widetilde\dist_{p_n,r})\]
are semiconcave and continuously differentiable in $A^\circ$.

Suppose that a sequence of functions $\phi_i$ is $C^2$-converging to $\phi$.
Choose $x=(x_1,\dots,x_n)$ in the domain of definition of $\phi$.
Note that $\phi_n$ and its partial derivatives up to order 2 are bounded;
fix a bound $c$.
Then in a neighborhood of $(x_1,\dots,x_n)$ we may choose $\psi_n$ that is uniquely defined by $\psi_n(x)=0$, $\partial_i\psi_n(x)=2\cdot c$, $\partial_i\partial_j\psi_n\equiv 0$ for $i\ne j$, and $\partial_i^2\psi_n\equiv -d$ for a large constant $d$. 
In this case, $\chi_n=\psi_n-\phi_n$ is concave.
Moreover, $C^2$-convergence of $\phi_n$ implies convergence of $\psi_n$ and $\chi_n$.
Hence, the second statement follows.
\qeds

The definition of DC-convergence extends naturally to sequences of functions defined on a fixed domain $\Omega\subset \RR^m$.
The proof of the following statement is a straightforward modification of \cite[Section 3]{PerDC}:

\begin{thm}{Proposition}\label{prop:DC-conv}
Let $M_n\smooths{} A$;
choose a common chart $\bm{x}_n\:U_n\zz\subset M_n\to \Omega$, $\bm{x}\:U\subset A\zz\to \Omega\subset \RR^m$.

Consider functions $f_n$ and $f$ defined on $U_n$ and $U$ respectively.
Then $f_n\DCto f$ if and only if $f_n\circ \bm{x}_n^{-1}\DCto f\circ \bm{x}^{-1}$
\end{thm}

The following statement follows from the lemma in \cite[Section 4]{PerDC}.

\begin{thm}{Proposition}\label{metricBV+}
Let $A$ be an $m$-dimensional Alexandrov space and $\bm{x}\:U\to \RR^m$ --- a smoothed chart for $U\subset A$.
Denote by $g_{ij}$ components of metric tensors in this chart
and by $g^{ij}$ components of the inverse matrix. 
Let $f\:U\to \RR$ be a DC$_0$ function.

Then the partial derivatives $\partial_k g_{ij}$, $\partial_k g^{ij}$, $\partial_i\partial_jf$ are Radon measures on $A$ that vanish on $\bm{x}^{-1}(A'\cup A'')$.
\end{thm}

\begin{thm}{Theorem}\label{metricBV}
Let $M_n\smooths{} A$, $\dim A=m$; choose a common chart $\bm{x}_n$ defined on $U_n\subset M_n$, $\bm{x}$ defined on $U\subset A$ with a common range $\Omega\subset \RR^m$.
Denote by $g_{ij,n}$ components of metric tensors in this chart
and by $g^{ij}_n$ components of the inverse matrix. 
Let $f_n\:U_n\to \RR$ be a sequence of DC function that DC-converges to a DC$_0$ function $f\:U\to \RR$.
Then partial derivatives $\partial_k g_{ij,n}$, $\partial_k g^{ij}_n$, $\partial_i\partial_jf_n$ weakly converge to the Radon measures $\partial_k g_{ij}$, $\partial_k g^{ij}$, $\partial_i\partial_jf$ described in \ref{metricBV+}.
\end{thm}

By \ref{obs:test-DC}, the theorem applies to any test-converging sequence $f_n\zz\testto f$.
In the proof, we modify the argument in \cite[Section 4]{PerDC} slightly.

\parit{Proof.}
Let's start with the partial derivatives of metric tensors.
In \cite[Subsection 4.2]{PerDC}, it was shown that components of metric tensors can be expressed as a rational function of partial derivatives of distance functions to a finite collection of points.
The distance functions are semiconcave, in particular DC.

The base points $p_{i,n}\in M_n$ of these distance functions can be chosen so that they converge to some point $p_i\in A$.
In this case, the distance functions are DC-converging.
Now, applying \ref{prop:DC-conv}, we get the statement.

The case of $\partial_i\partial_jf_n$ is similar.
\qeds

\section{Bi-Lipschitz covering}\label{sec:bilip}

In this section we will prove Lemma~\ref{lem:A-prime-Q}.
A more general version of the lemma can be proved along the same lines as Lemma~11.1 in \cite{simon}.

Note that the lemma follows from the next proposition.

\begin{thm}{Proposition}\label{prop:Q-covering}
Let $A$ be an $m$-dimensional Alexandrov space with curvature at least $-1$ and $p\in A'$.
Then there is a compact set $Q$ such that 
\begin{enumerate}[(i)]
 \item $Q$ admits a bi-Lipschitz embedding into $\RR^{m-2}$ and
 \item there is a neighborhood $U\ni p$ and $\eps>0$ such that $q\in Q$ for any point $q\in U\cap A'$ such that 
 \[\theta(q)<\theta(p)+\eps.\]
\end{enumerate}
\end{thm}

Let $x$ be a point in an Alexandrov space $A$ with curvature at least $-1$.
Recall that Bishop--Gromov inequality implies that 
\[\frac{\vol^m B(x,R)_A}{\vol^m B(\tilde x,R)_{\HH^m}}
\le
\frac{\vol^{m-1} \Sigma_x}{\vol^{m-1} \SS^{m-1}}\]
for any $R>0$; here $\HH^m$ denotes the $m$-dimensional hyperbolic space.
The following lemma makes this inequality more precise. 

\begin{thm}{Lemma}
Let $x$ be a point in an $m$-dimensional Alexandrov space $A$ with curvature at least $-1$.
Suppose $y\in A$ is a point such that $|x-y|<R$ and $\measuredangle\hinge yxz<\pi-\eps$ for any point $z$.
Then
\[\frac{\vol^m B(x,R)_A}{\vol^m B(R)_{\HH^m}}
\le
(1-\delta)\cdot\frac{\vol^{m-1} \Sigma_x}{\vol^{m-1} \SS^{m-1}},\]
where $\delta$ is a positive number that depends on $m$, $|x-y|$, $R$ and $\eps$.
 
\end{thm}

\parit{Proof.}
To simplify the presentation we will assume that $A$ is nonnegatively curved;
it is straightforward to adapt the proof to the general case.
In this case, we need to show that 
\[\frac{\vol^m B(x,R)_A}{\vol^m B(R)_{\RR^m}}
\le 
(1-\delta)\cdot\frac{\vol^{m-1} \Sigma_x}{\vol^{m-1} \SS^{m-1}},\]

Let us denote by $\tilde p$ a vector in $\T_x$ that is tangent to a geodesic path $\gamma\:[0,1]\zz\to A$ from $x$ to~$p$.
By comparison, the map $p\mapsto \tilde p$ is a distance-noncontracting map.

Since $\measuredangle\hinge yxz<\pi-\eps$ for any $z$, the image of the map $p\mapsto \tilde p$ does not include points in a cone $C$ behind $\tilde y$ of angle $\eps$.
It follows that 
\[\vol^m(B(0,R)_{\T_x}\setminus C)\ge \vol^m(B(x,R)_A)\]
for any $R>0$.

\begin{wrapfigure}{r}{35 mm}
\vskip-0mm
\centering
\includegraphics{mppics/pic-200}
\vskip0mm
\end{wrapfigure}

Since $R>|x-y|$, the intersection $C\zz\cap B(0,R)_{\T_x}$ includes a ball of a certain radius $r>0$ that can be found in terms of $|x-y|$, $R$ and $\eps$.
By Bishop--Gromov inequality, we get $\delta=\delta( m, |x-y|, R, \eps)>0$ such that 
\[\frac{\vol^m(C\cap B(0,R)_{\T_x})}{\vol^m(B(0,R)_{\T_x})}>\delta.\]

Further, observe that 
\[\frac{\vol^m(B(0,R)_{\T_x})}{\vol^m(B(0,R)_{\RR^m})}
=
\frac{\vol^{m-1} \Sigma_x}{\vol^{m-1} \SS^{m-1}}.
\]
--- whence the lemma.
\qeds

\parit{Proof of \ref{prop:Q-covering}.}
Since the tangent cone at $p$ has $\RR^{m-2}$-factor, 
we can choose points $a_1,\zz\dots,a_{m-2},b_1,\zz\dots,b_{m-2}$ that are $\delta$-strainers of $p$ for arbitrary $\delta>0$.
The corresponding distance map $s\:x\zz\mapsto (|a_1-x|,\zz\dots,|x-a_{m-2}|)$ is an almost submersion of a neighborhood $U\ni p$ to $\RR^{m-2}$.
Choose small $\eps>0$ and set 
\[Q'=\set{x\in U\cap A'}{\theta(x)<\theta(p)+\eps}.\]
Let us show that $s|_{Q'}$ is bi-Lipschitz.
Once it is done, passing to the closure $Q=\bar Q'$ gives the required set. 

Note that for some $R>0$ the ball $B(p,10\cdot R)_A$ is almost isometric to the ball $B(0,10\cdot R)_{\T_p}$
and we can assume that $U\subset B(p,R)_A$.
By the volume convergence (see \cite[10.8]{BGP}) and Bishop--Gromov inequality, we can assume that 
\[\vol^m B(x,R)_A>\tfrac{\theta(p)-\eps}{2\cdot \pi}\cdot\vol^m B(0,R)_{\HH^m}\]
for any $x\in U$; here $\HH^m$ denotes the $m$-dimensional hyperbolic space. 

Assume $x$ and $y$ in $Q'$.
Since $\eps$ is small, the lemma implies that there is $z\in A$ such that $\measuredangle\hinge yxz$ is near $\pi$.
It follows that $\dir yx$ lies very close to the $\RR^{m-2}$-factor in~$\T_y$.
The same way we can show that $\dir xy$ lies very close to the $\RR^{m-2}$-factor in~$\T_x$.
In other words $[xy]$ lies nearly horizontally with respect to almost submetry $s$.
In particular,  
\[|s(x)-s(y)|_{\RR^{m-2}} \lessgtr \lambda^{\pm1}\cdot |x-y|_{A}\]
some constant $\lambda>1$.
(In fact, we can take $\lambda$ arbitrarily close to $1$, but we do not need it.)
\qeds

\section{Localization}\label{sec:local}

In this section we formulate a local version of the main theorem.
This version is more general, but its proof requires just a slight change of language.
A couple of times we had to use this local version in the proof.
In a perfect world, we had to rewire the whole paper using this language.
However, this is not a principle moment,
so we decided to keep the paper more readable at the cost of being not fully rigorous.
A more systematic discussion of this topic is given in \cite{LNep}.

First, we need to define Alexandrov region;
its main example is an open set in Alexandrov space.

\begin{thm}{Definition}
Let $A$ be a locally compact metric space. 
We say that a point  $p\in A$
is \emph{$\eps$-inner} if
the closed ball $\bar B(x,2\cdot\eps)$ is compact.
\end{thm}

\begin{thm}{Definition}
We say that a locally compact inner metric space $A$ of finite Hausdorff dimension is an \emph{Alexandrov region} if any point has a neighborhood where the Alexandrov comparison for curvature $\ge -1$ holds.

The \emph{comparison radius} $r_c(p)$ for $p\in A$ is defined as the maximal number $r$ such that $p$ is $r$-inner point and Alexandrov comparison for curvature $\ge -1$ holds in $B(x,r)$.
\end{thm}

Any point $p$ in an Alexandrov region admits a convex neighborhood.
Moreover, its size can be controlled in terms of dimension, $r_c(p)$, and a lower bound on the volume of ball $B(p,r_c)$.
The construction is the same as for Alexandrov space \cite[4.3]{perelman-petrunin}.

By the globalization theorem (see, for example, \cite{AKP}), a compact convex subset in an Alexandrov region is an Alexandrov space.
So the statement above makes it possible to apply most of the arguments and constructions for Alexandrov spaces to Alexandrov regions. 
Moreover, in the case when an Alexandrov region is a Riemannian manifold (possibly noncomplete) it is possible to take the doubling of a convex neighborhood from the proposition and smooth it with almost the same lower curvature bound.
This allows us to apply the main result from \cite{petrunin-SC}, where the complete manifold can be replaced by a convex domain in a possibly open manifold. 

Further, let us define a local version of smoothing.
Let us denote by
$\M_{\ge -1}^m$ a class of $m$-dimensional Riemannian 
manifolds without boundary, but possibly non-complete, with sectional curvature bounded
from below by $-1$.

\begin{thm}{Definition}
Let 
$M_n\in\M_{\ge -1}^m$ (with corresponding intrinsic metric)
converge in Gromov--Hausdorff sense to some metric space $A$ via
approximation.
Suppose that $M_n\ni x_n\to x\in A$, $\dim A=m$,
and $r_c(x_n)\ge R>0$.
Let 
$U_n=B(x_n,R)_{M_n}$.
Then we say that $U_n$ is a local smoothing of $U=B(x,R)_A$ (briefly, $U_n\smooths{} U$).
\end{thm}

It is straightforward to redefine test functions and weak convergence for local smoothings.
Using this language we can make a local version for each statement in this paper, the proofs go without changes.
As a result, we get the following local version of the main theorem \ref{main}.
 
\begin{thm}{Local version of the main theorem}\label{mainloc}
Let   
$M_n\in\M_{\ge -1}^m$,
$M_n\GHto A$, $U_n\subset M_m$,
  $U\subset A$, and $U_n\smooths{} U$ be a local smoothing.
  
Denote by $\qm_n$ the dual measure-valued curvature tensor on $U_n$.
Then there is a measure-valued tensor $\qm$ on $U$ such that $\qm_n\rightharpoonup \qm$.
\end{thm}

{\sloppy

\def\emph{\textit}
\printbibliography[heading=bibintoc]
\fussy
}

\Addresses

\end{document}